\documentclass[11pt,a4paper,leqno]{amsart}

\usepackage[colorlinks,linkcolor=blue,citecolor=blue,urlcolor=teal, linktocpage=true, breaklinks=true]{hyperref}

\usepackage{amsmath,amsfonts,amssymb,amsthm}
\usepackage[english]{babel}
\usepackage{enumerate}

\usepackage{lipsum}

\newtheorem{theor}{Theorem}[section]
\newtheorem{coro}[theor]{Corollary}
\newtheorem{prop}[theor]{Proposition}
\newtheorem{lemma}[theor]{Lemma}

\theoremstyle{definition}
\newtheorem{defin}[theor]{Definition}
\newtheorem{remark}[theor]{Remark}

\numberwithin{equation}{section}

\usepackage[dvipsnames]{xcolor}

\def\R{\mathbb{R}}
\def\C{\mathbb{C}}
\def\D{\mathbb{D}}
\def\Z{\mathbb{Z}}

\def\E{\mathbf{E}}
\def\V{\mathbf{V}}
\def\P{\mathbf{P}}

\def\K{\mathcal{K }}

\newcommand{\uno}{\mathop{\textbf{1}}}

\title[Moments and growth of power series distributions]{Growth of power series with nonnegative coefficients, and moments of power series distributions}

\author[A. Cant\'on]{Alicia Cant\'on}
\address[Alicia Cant\'on]{Departamento de Matem\'atica e
Inform\'atica Aplicadas a las Ingenier\'{\i}as Civil y Naval,
Universidad Polit\'ecnica de Madrid, Spain}
\email{alicia.canton@upm.es}

\author[J.\,L. Fern\'andez]{Jos\'e L. Fern\'andez}
\address[Jos\'e L. Fern\'andez]{Departamento de Matem\'aticas, Universidad Aut\'onoma de Madrid, Spain.}
\email{joseluis.fernandez@uam.es, jlfernandez@akusmatika.org}

\author[P. Fern\'andez]{Pablo Fern\'andez}
\address[Pablo Fern\'andez]{Departamento de Matem\'aticas, Universidad Aut\'onoma de Madrid, Spain}
\email{pablo.fernandez@uam.es}

\author[V.\,J. Maci\'a]{V\'{\i}ctor J. Maci\'a}
\address[V\'{\i}ctor J. Maci\'a]{Departamento de Matem\'aticas, Universidad Aut\'onoma de Madrid, Spain.}
\email{victor.macia@uam.es}


\begin{document}

\begin{abstract}
Any power series with nonnegative coefficients has an associated family of probability distributions supported on the
nonnegative integers. There is a close connection between the function theoretic properties of the power series and
the moments of the family of distributions. In this paper, we describe that interplay, provide simpler proofs of some
known results by emphasizing the probabilistic perspective, and present some new theorems.
\end{abstract}

\maketitle

\begin{flushright}\emph{Dedicated to the memory of Luis B\'aez-Duarte.}
\end{flushright}



\setcounter{tocdepth}{2}


\section{Introduction}

Every power series $f(z)=\sum_{n=0}^\infty a_n z^n$, with nonnegative coefficients and
radius of convergence $R>0$, generates a family $(X_t)_{t\in[0,R)}$
of probability distributions supported on the nonnegative integers by defining random
variables $X_t$ via
\[
X_0\equiv 0,\qquad \P(X_t=n)=\frac{a_n t^n}{f(t)}, \quad \text{for $t\in (0,R)$ and $n \ge 0$}.
\]
The systematic study of these \emph{power series distributions} appears to have originated with Kosambi \cite{Kosambi1}
in 1949 and Noak \cite{Noak} in 1950, as a unifying model for several discrete distributions arising in Statistics.
See Section~2.2 in~\cite{JohnsonUnivariate} and references therein.

We will restrict our study to the class $\mathcal{K}$ of \emph{nonconstant} power series with nonnegative coefficients,
positive radius of convergence, \emph{and with $a_0>0$}. With this, we avoid the inconvenient case where $f$ is a monomial,
or a constant.

We will refer to the family of random variables $(X_t)_{t \in [0,R)}$ associated to~$f$ in $\mathcal{K}$ as the
\emph{Khinchin family} of $f$, a term and a whole framework which originates in the work of Rosenbloom~\cite{Rosenbloom}, see also \cite{Schumitzky}.

This paper deals with the beautiful interplay between function theoretical properties of the holomorphic function given by $f$ and
certain probabilistic properties of the family $(X_t)$ as~$t\uparrow R$, in particular, the behaviour as $t\uparrow R$ of
the mean $m_f(t)=\E(X_t)$, of the variance $\sigma_f^2(t)=\V(X_t)$ or, in general, of the moments $\E(X_t^p)$, for $p>0$.

For instance, in Section \ref{section:growth of variance} we will provide simple proofs, \emph{of probabilistic nature},
of some lower bounds and asymptotic lower bounds of $\sigma_f(t)$, as $t\uparrow R$, quantifying Hadamard's three lines theorem, due
to Hayman~\cite{Hayman2}, {Bo\u{\i}chuk} and Gol'dberg~\cite{BoichukGoldberg},  Abi-Khuzzam~\cite{Abi2} and others. Later, see
Theorem~\ref{teor:orden y m} and Proposition~\ref{prop:orden y sigma cuad partod por m}, we will show, for entire functions~$f$
in $\mathcal{K}$, how the growth of the mean $m_f(t)$ and, in general, of the moments $\E(X_t^p)$, with $p>0$, and of the
quotient $\sigma^2_f(t)/m_f(t)$, relate to the order of $f$, generalizing some results of P\'olya and Szeg\H{o}, and of B\'aez-Duarte.

We say that the  family $(X_t)_{t \in [0,R)}$ is a \emph{clan} if  $\lim_{t\uparrow R} \sigma_f(t)/m_f(t)=0$ or, equivalently,
if $\lim_{t\uparrow R} \E(X_t^2)/\E(X_t)^2=1$. For clans, $X_t$ concentrates around its mean as $t\uparrow R$, in the sense that
the normalized variable $X_t/\E(X_t)$ tends to the constant 1 in probability. As a matter of fact, being a clan was already
considered, not with this name, by Hayman in~\cite{Hayman} as a property enjoyed by what are nowadays called Hayman (admissible)
functions; clans also appear at least in P{\' o}lya--Szeg\H{o} (see items~70 and~71 in~\cite{PolyaSzego}), and in work of
Simi{\'c}~\cite{Simic,Simic0}.

Section \ref{section:clans} contains a detailed study of clans from a function theoretical point of view.
It turns out, Theorem \ref{teor:asintotico de momentos si clan}, that for clans, $\lim_{t\uparrow R} \E(X_t^p)/\E(X_t)^p=1$,
for any $p>0$. For instance, the partition function $P(z)=\prod_{k=1}^\infty 1/(1{-}z^k)$ is a clan, and for the family
$(X_t)_{t>0}$ associated to $P$ we obtain readily, for any $p >0$, that
$\E(X_t^p)\sim \zeta(2)^p/(1-t)^{2p}$, as $t \uparrow 1$.

For entire functions in $\K$, P\'olya and Szeg{\H{o}} \cite{PolyaSzego} showed that entire functions with nonnegative coefficients
of finite order $\rho$ such that $\lim_{t \to \infty} \ln f(t)/t^\rho$ exists and is positive are clans; we show, more generally,
see Theorem \ref{theor:proximate order and clans} in Section \ref{section:entire functions}, that entire functions in~$\K$ of
regular growth (in a precise sense) are clans.

The entire gap series with nonnegative coefficients presented in Section~\ref{section:entire gap series} furnish the basic
examples of entire functions of any order $\rho$, with $0 \le \rho \le +\infty$, which are not clans. A classical result of Pfluger
and P\'olya~\cite{PflugerPolya} ensures that these entire gap series have no Borel exceptional values.
We show in Theorem \ref{teor:Borel exceptional and clans} that entire functions in $\mathcal{K}$ with one Borel exceptional value
are always clans.

Finally, some open questions are discussed in Section \ref{section:questions}.

\medskip

\noindent\textbf{Notation}. Generically, along this paper, we use $\E(Z)$ and $\V(Z)$ to denote expectation and variance
of a random variable $Z$, and $\P(A)$ to denote probability of the event $A$ (in the appropriate probability space).

For positive functions $f(t)$ and $g(t)$, defined in an interval $[0,R)$, with $0<R\le +\infty$, the notation $f(t)\sim g(t)$
as $t\uparrow R$ means that $\lim_{t\uparrow R} f(t)/g(t) = 1$, while $f(t)=o(g(t))$ as $t\uparrow R$ means that
$\lim_{t\uparrow R} f(t)/g(t) = 0$. We write $f(t)=O(g(t))$ as $t\uparrow R$ if $f(t)\le C g(t)$ for a certain constant $C>0$ and
for $t$ close enough to $R$; and $f(t)\asymp g(t)$ as $t\uparrow R$ if $cg(t)\le f(t)\le C g(t)$ for certain constants $c,C>0$
and for $t$ close enough to $R$.

\section{Khinchin families and power series distributions}\label{section:Khinchin families}

Keeping the notation and definitions {of} 
 \cite{CFFM1}, we denote by $\K$ the class of nonconstant power series
\begin{equation}
\label{eq: serie de potencias}
f(z)=\sum_{n=0}^\infty a_n z^n
\end{equation}
with positive radius of convergence $R\le \infty$, which have \emph{nonnegative} Taylor coefficients, and with $a_0>0$. Thus,
one coefficient of the power series $f$ other that $a_0$ is also nonzero. Notice that $f(t)>0$ for each $t \in[0,R)$.
We shall resort occasionally to the fact that for~$f$ in $\K$, we have that
\begin{equation}
\label{eq:max in t}
\max\{|f(z)|: |z|\le t\}=\max\{|f(z)|: |z|= t\}=f(t), \quad \text{for every $t \in [0,R)$}.
\end{equation}

To any power series $f$ in $\K$, we may associate a whole family of probability distributions supported on the
nonnegative integers $\{0,1, \ldots\}$, indexed in the interval $[0,R)$, by specifying a family of random variables
$(X_t)_{t \in[0,R)}$ with values in $\{0, 1, \ldots\}$ and with probability mass functions given by
\[
\P(X_t=n)=\frac{a_n t^n}{f(t)}, \quad \text{for each $n \ge 0$ and $t\in (0,R)$},
\]
and with $X_0\equiv 0$. Since $f$ has at least two coefficients which are not zero, each variable~$X_t$, for
$t \in (0,R)$, is a nonconstant random variable.

The family $(X_t)_{t\in[0,R)}$ is the \emph{family of probability distributions associated to the power series $f$}.
For general background on power series probability distributions, we refer the reader to~\cite{JohnsonUnivariate},
Section~2.2.

The terminology of \emph{Khinchin families} for $(X_t)_{t\in[0,R)}$ is used in this paper to emphasize that the focus
is placed on the behaviour as $t \uparrow R$ of the variables~$(X_t)$ and of the function $f(t)$, and in the connection
between the probabilistic and the function theoretical aspects.

The framework of Khinchin families arises in work of Hayman \cite{Hayman}, Rosenbloom~\cite{Rosenbloom},
B\'aez-Duarte~\cite{BaezDuarte} and others. It has been pursued recently by the authors in~\cite{CFFM1}
and~\cite{CFFM2}. The main aim of that framework has been to study, in a unified manner,
asymptotic formulas for the coefficients of generating power series, like the Hardy--Ramanujan asymptotic
formula for the number of partitions of numbers with $f(z)=\prod_{n=1}^\infty 1/(1-z^n)$, the Moser--Wyman
asymptotic formula for the number of partitions of sets with $f(z)=e^{e^z-1}$, and many others.

The presentation of this paper is independent of the papers above. For the general theory of Khinchin families, we refer
the reader to \cite{CFFM1} (and also to \cite{CFFM2}). The summary of a few relevant notions follows.

\subsection{Basic families}
The most basic examples of Khinchin families $(X_t)_{t\in [0,R)}$ associated to power series $f\in \K$ are
\begin{itemize}
\item the \emph{Bernoulli family}, associated to $f(z)=1+z$, where for each $t >0$,
the variable~$X_t$ is a Bernoulli variable with parameter $p=t/(1+t)$;

\item the \emph{geometric family}, associated to $f(z)=1/(1-z)$, where for each $t \in(0,1)$,
the variable $X_t$ follows a geometric distribution with parameter $p=1-t$;

\item for integer $N \ge 1$, the \emph{binomial family}, associated to $f(z)=(1+z)^N$, where
for each $t >0$, the variable $X_t$ follows a binomial distribution 
 of parameters $N$ and $p=t/(1+t)$;

\item for integer $N \ge 1$, the \emph{negative binomial family}, associated to $f(z)=1/(1-z)^N$,
where for each $t \in (0,1)$, the variable $X_t$ follows a negative binomial distribution of
parameters $N$ and $p=1-t$;

\item and the \emph{Poisson family}, associated to $e^z$, where for each $t >0$, the variable $X_t$
follows a Poisson distribution of parameter $t$.
\end{itemize}

We mention also the families associated to
\begin{itemize}
\item the (ordinary) generating function of partitions of integers:
$\prod_{n=1}^\infty 1/(1-z^n)$, for $|z|<1$, 

\item and the (exponential) generating function of the Bell numbers (number
of partitions of sets): $e^{e^z-1}$, which is an entire power series.
\end{itemize}

\subsection{Moments}\label{section:comparison of moments of Khinchin families}

Let $f $ be a power series in $\K$, with radius of convergence $R>0$, and
let $(X_t)_{t\in[0,R)}$ be its associated Khinchin family. For $t \in [0,R)$,
the mean and variance of $X_t$ will be denoted by $m_f(t)=\E(X_t)$ and $\sigma_f^2(t)=\V(X_t)$,
respectively. In terms of $f$ itself, we may write
\begin{equation}
\label{eq:m and sigma in terms of f}
m_f(t)=\frac{t f^\prime(t)}{f(t)}=t \,(\ln f)^\prime(t) \quad \text{and} \quad \sigma_f^2(t)=t m_f^\prime(t),
\quad \text{for $t \in [0,R)$}.
\end{equation}

More generally, for $\beta \ge 0$ and $t\in[0,R)$, \emph{the moment $\E(X_t^\beta)$ of exponent $\beta$ of\/ $X_t$}
may be written in terms of the power series  $f$ as
\[
\E(X_t^\beta)=\sum_{n=0}^\infty n^\beta \,\frac{a_n t^n}{f(t)}\cdot
\]
We shall be particularly interested in the comparison of the moment $\E(X_t^\beta)$ of exponent~$\beta$
with the mean $m_f(t)=\E(X_t)$, which acts as a sort of unit, as $t\uparrow R$.

Occasionally, we shall resort to \emph{factorial moments} of the $X_t$. For integer $k \ge 0$, we
denote $x^{\underline{k}}:=x(x-1)\cdots(x-k+1)$, the $k$-th falling factorial of the number $x$.
For $k=0$, it is agreed that $x^{\underline{0}}=1$. The $k$-th factorial moment of $X_t$ is defined as
\[
\E(X_t^{\underline{k}})=\E(X_t\, (X_t-1)\cdots (X_t-k+1)),
\]
with the convention that $\E(X_t^{\underline{0}})=1$.
Observe that the random variable $X_t^{\underline{k}}$ given
by $X_t^{\underline{k}}:=X_t\, (X_t-1)\cdots (X_t-k+1)$ takes
values in $\{0, 1, \ldots\}$, and that $\E(X_t^{\underline{k}})\ge 0$,
for any~$k\ge 0$.

In terms of the power series $f$ itself, the $k$-th factorial moment is given by
\begin{equation}
\label{eq:formula factorial moments}
\E(X_t^{\underline{k}})
=
\frac{1}{f(t)}\sum_{n=k}^\infty n(n-1)\cdots (n-k+1) \,a_n\, t^n=\frac{t^k f^{(k)}(t)}{f(t)},
\quad \text{for $t \in [0,R)$.}
\end{equation}

Stirling numbers of the second kind $S(n,j)$ relate powers with falling factorials of a number $x$:
\begin{equation*}
x^k=\sum_{j=0}^k S(k,j) \, x^{\underline{j}}\,, \quad \text{for $k \ge 1$.}
\end{equation*}
This translates into the following relation between moments and factorial moments, see, for instance,
equation~(11) of~\cite{Pitman}:
\begin{equation}
\label{eq:moments factorial stirling}
\E(X_t^k)=\sum_{j=0}^k S(k,j) \, \E(X_t^{\underline{j}}),
\quad \text{for any integer $k \ge 1$ and $t \in (0,R)$.}
\end{equation}
Recall that $S(k,k)=1$, for any $k \ge 1$.

 We describe next the moments and asymptotics of the moments of the most basic families.

\subsubsection{Moments of geometric and negative binomial variables}

\begin{prop}\label{prop:moments geometric}
If\/ $N \ge 1$ and $(X_t)_{t\in[0,1)}$ is the Khinchin family of\/ $f(z)=1/(1-z)^N$, then
\[
m_f(t)=N\frac{t}{1-t} \quad\text{{and}}\quad  \sigma_f^2(t)=N\frac{t}{(1-t)^2}\,, \quad \text{for $t \in [0,1)$,}
\]
and for any $\beta\ge 0$, we have that
\[
\lim_{t\uparrow 1}\dfrac{\E(X_t^\beta)}{\E(X_t)^\beta}=\frac{\Gamma(\beta + N)}{\Gamma(N) N^\beta}\cdot
\]
\end{prop}

This comparison of moments is a direct consequence of the following asymptotic formula (valid for any $\beta>0$):
\begin{equation*}
\sum_{n=1}^\infty n^{\beta-1} t^n \sim \Gamma(\beta)\,\frac{1}{(1-t)^{\beta}}, \quad \text{as $t\uparrow 1$},
\end{equation*}
which in turn follows from the binomial expansion and Stirling's formula.

\subsubsection{Moments of Bernoulli and binomial variables}

For the Khinchin family $(X_t)_{t\ge 0}$ of a \emph{polynomial $f$ of degree $N$}, we have that $X_t$ tends
in distribution to the constant $N$ as $t \to \infty$, and also that for $\beta \ge 0$,
\[
\lim_{t\to \infty} \E(X_t^\beta)=N^\beta=\lim_{t\to \infty} \E(X_t)^\beta.
\]

Particular instances of this polynomial case are the Bernoulli 
 and the binomial families, which are associated, respectively, to $f(z)=1+z$ (of degree 1) and to
 $f(z)=(1+z)^N$ (of degree $N$).

\subsubsection{Moments of Poisson variables}
The Poisson family $(X_t)_{t \ge 0}$ is associated to the exponential function $e^z$. For $t>0$,
the variable $X_t$ is a Poisson variable with parameter~$t$, and thus $\E(X_t)=t$ and $\V(X_t)=t$.

For $\beta >0$ and $t >0$, the moment of exponent $\beta$ of $X_t$ is given by
\[
\E(X_t^\beta)=\sum_{n=0}^\infty n^\beta \,\frac{t^n e^{-t}}{n!}\cdot
\]

\begin{prop}
\label{prop:positive moments of poisson}
For the family $(X_t)_{t\ge 0}$ associated to $f(z)=e^z$, and for any $\beta>0$, we have
\[\E(X_t^\beta)\sim t^\beta=\E(X_t)^\beta, \quad \text{as $t \to \infty$}.\]
\end{prop}

An application of Jensen's inequality will allow to reduce the discussion below to the case of integer exponents.

\begin{lemma}
\label{lema:from integer moments to general moments}
Let $0<R\le\infty$, and let $(U_t)_{t\in[0,R)}$ be a family of nonnegative random variables
such that for some $p>1$,
\[
\E(U_t^p)\sim \E(U_t)^p, \quad\text{as } t\uparrow R.
\]
Then, for any $\beta\in (0,p]$,
\[
\E(U_t^\beta)\sim \E(U_t)^\beta, \quad\text{as } t\uparrow R.
\]
\end{lemma}

\begin{proof}
Let $1<\beta\le p$. We have from Jensen's inequality that
\[
 \E(U_t)\le \E(U_t^\beta)^{1/\beta}\le \E(U_t^p)^{1/p}, \quad \text{for any $t \in (0,R)$}.
\]
The result follows in this case, since, by hypothesis, $\E(U_t^p)^{1/p}\sim \E(U_t)$, as $t \uparrow R$.

For $\beta \in (0,1)$, Jensen's inequality gives that
\[
\E(U_t^\beta)\le \E(U_t)^\beta, \quad \text{for any $t \in (0,R)$}.
\]
Consider now $u\in (0,1)$ such that $1=\beta u +p (1-u)$.
By H\"older's inequality,
\[
 \E(U_t)\le \E(U_t^\beta)^u \, \E(U_t^p)^{1-u}, \quad \text{for any $t \in (0,R)$}.
\]
The result follows from these two inequalities, since, by hypothesis and the definition of~$u$, we have
that $\E(U_t^p)^{1-u}\sim \E(U_t)^{1-\beta u}$, as $t \uparrow R$.
 \end{proof}

\begin{proof}[Proof of Proposition \textup{\ref{prop:positive moments of poisson}}]
As $f(t)=e^t$, it follows from \eqref{eq:formula factorial moments} that the factorial moments of $X_t$
are given by $\E(X_t^{\underline{j}})=t^j$, for any integer $j\ge 0$ and any $t >0$.
On account of~\eqref{eq:moments factorial stirling}, we have that, for each integer $k \ge 1$,
\[
\E(X_t^k)=\sum_{j=0}^k S(k,j) \,t^j\,.
\]
Since $S(k,k)=1$, for $k\ge 1$, we obtain that for each integer $k \ge 1$,
\[
\E(X_t^k)\sim t^k, \quad \text{as $t \to \infty$}.
\]
Lemma \ref{lema:from integer moments to general moments} finishes the proof.
\end{proof}

\subsubsection{Mean and variance of the partition function and of the Bell function}
Beyond the most basic families, consider now the partition function $P(z)$ given by the infinite product
\[
P(z)=\prod_{k=1}^\infty \frac{1}{1-z^k}=\sum_{n=0}^\infty p(n)\, z^n , \quad \text{for $|z|\le 1$.}
\]
The coefficient $p(n)$ is the number of partitions of the integer $n$.

The mean $m_P(t)$ and the variance $\sigma_P^2(t)$ are given, on account of \eqref{eq:m and sigma in terms of f}, by
\[
m_P(t)=\sum_{n=1}^\infty \frac{n t^n}{1-t^n} \quad\text{and}\quad  \sigma_P^2(t)=\sum_{n=1}^\infty \frac{n^2 t^n}{(1-t^n)^2},
\quad \text{for $t\in[0,1)$.}
\]
By Euler summation, one may obtain convenient asymptotic formulas describing their behaviour as $t\uparrow 1$:
\begin{equation}
\label{eq:mean and variance partition}
m_P(t)\sim \frac{\zeta(2)}{(1-t)^2}\quad\text{and}\quad \sigma_P^2(t)\sim \frac{2\zeta(2)}{(1-t)^3}, \quad \text{as $t \uparrow 1$.}
\end{equation}
See, for instance, Section 6.1 of \cite{CFFM1}.

We will see later on, see Section \ref{seccion:examples of clans}, that $P$ is a clan.
Thus Theorem~\ref{teor:asintotico de momentos si clan} would claim, for the family $(X_t)$ associated to the partition
function $P$, that for any $\beta >0$ we have that $\E(X_t^\beta)\sim \E(X_t)^\beta$, as $t \uparrow 1$,
and therefore that
\[
{\E(X_t^\beta)\sim\frac{\zeta(2)^\beta}{(1-t)^{2\beta}}} , \quad \text{as $t \uparrow 1$,}
\]
for any $\beta >0$.

The exponential generating function $B(z)$ of the Bell numbers $B_n$ is given by
\[
B(z)=e^{e^z-1}=\sum_{n=0}^\infty B(n)\frac{z^n}{n!} , \quad \text{for $z \in \C$.}
\]
The coefficient $B(n)$ is the number of partitions of the set $\{1,\ldots, n\}$. In this case,
the mean and variance admit simple formulas:
\begin{equation}
\label{eq:mean and variance Bell}
m_B(t)=te^t\quad\text{and}\quad \sigma_B^2(t)=t(t+1) e^t , \quad \text{for $t >0$.}
\end{equation}

The Bell function is also a clan: for the moments of the family $(X_t)$ associated to the Bell function we have,
for any $\beta>0$, that
\[
\E(X_t^\beta)\sim t^\beta e^{\beta t} , \quad \text{as $t \to \infty$.}
\]

\subsubsection{Mean and variance of some canonical products}\label{section:mean variance canonical products}

Let  $(b_k)_{k \ge 1}$ be a sequence of positive numbers increasing to $\infty$ in such a way that $\sum_{k=1}^\infty 1/b_k<+\infty$. The \emph{canonical product} $f$ given by
\begin{equation}\label{eq:formula of canonical product}
f(z)=\prod_{k=1}^\infty \Big(1+\frac{z}{b_k}\Big), \quad \text{for $z \in \C$.}
\end{equation}
is an entire function in $\K$, {whose set of zeros is}
 $\{-b_k, k \ge 1\}$. {By Hadamard's factorization Theorem, these are \emph{all} the entire functions of genus 0 with only negative (real) zeros, normalized so that $f(0)=1$. See, for instance, Chapter 4 of~\cite{Boas}.

For the mean and variance functions of $f$, we have, from \eqref{eq:m and sigma in terms of f},  that
\begin{equation}
\label{eq:mean and variance of canonical products}
m_f(t)=\sum_{k=1}^\infty \frac{t}{t+b_k}\quad \text{and} \quad \sigma_f^2(t)=\sum_{k=1}^\infty \frac{b_k t}{(t+b_k)^2} ,
\quad \text{for $t \ge 0$.}
\end{equation}

These canonical products provide interesting examples of the behaviour of $\sigma_f^2(t)$, as we discuss now, and they will be used in forthcoming arguments (see {Sections~\ref{subsec:entire}~and~\ref{section:genus 0}}).
Denote with
\[
N(t)=\#\{k\ge 1: b_k \le t\},\quad \mbox{ for $t >0$,}
\]
the \emph{counting function of zeros} of $f$.
Thus $N(t)$ counts the number of zeros of $f$ in the disk $\D(0,t)$.

\smallskip

{The mean and variance functions of $f$ are readily comparable.}
\begin{lemma}\label{lema:variance less than mean of canonical products} For functions $f$ as in \eqref{eq:formula of canonical product}, we have
\begin{equation}
\label{eq:variance less than mean of canonical products}
\sigma_f^2(t)< m_f(t)<2\sigma_f^2(t)+N(t) , \quad \text{for any $t >0$.}
\end{equation}
\end{lemma}
\begin{proof}
On the one hand, observe that
\[
\frac{b_k t}{(t+b_k)^2}<\frac{t}{t+b_k}, \quad \mbox{for any $t >0$ and any $k \ge 1$.}
\]
On the other hand,
\[
\frac{t}{t+b_k}<1 \quad\text{for all $t>0$, but also}\quad
\frac{t}{t+b_k} <2 \,\dfrac{b_k t}{(t+b_k)^2}  \quad \text{if $b_k> t >0$.}
\]
The statement follows from these estimates and the formulas in \eqref{eq:mean and variance of canonical products}.
\end{proof}


{We follow the lead of Hayman in Theorem 4 of~\cite{Hayman2}, and show next how $\sigma_f$ of these canonical products depend upon the log spacing of the $b_k$.}

\begin{lemma}\label{lema:condition on bk}Consider a canonical product as in \eqref{eq:formula of canonical product}, and assume further that the sequence $(b_k)_{k \ge 1}$ satisfies
 \begin{equation}\label{eq:condition on bk} b_{k+1}\ge 2 b_k, \quad \mbox{for any $k \ge 1$}.
 \end{equation}
 For each $n \ge 2$, we let $I_n$ denote the interval
\begin{equation}\label{eq:def de In}
I_n:=[\sqrt{b_{n-1}b_n}, \sqrt{b_n b_{n+1}}].
\end{equation}
Then, for any $n \ge 2$,
\begin{equation}\label{eq:bounds sigma productos canonicos uniform}
\begin{aligned}\sup_{t\in I_n}\sigma_f^2(t)&\le \frac{1}{4} +4 \max\big\{\sqrt{{b_n}/{b_{n+1}}},\,\sqrt{{b_{n-1}}/{b_{n}}}\big\},
\\
\inf_{t\in I_n}\sigma_f^2(t)&\ge \frac{1}{4}\min\big\{\sqrt{{b_n}/{b_{n+1}}},\,\sqrt{{b_{n-1}}/{b_{n}}}\big\}.
\end{aligned}
\end{equation}
\end{lemma}
\begin{proof}
Consider the positive function
\begin{equation}\label{eq:def de varphi}
\varphi(x)=\frac{x}{(1+x)^2}, \quad\text{for $x >0$,}
\end{equation}
 which attains a maximum value of $1/4$ at $x=1$, and satisfies
\begin{equation}
\label{eq:estimates phi}
\frac{x}{4}<\varphi(x)<x, \quad \mbox{for $x \in (0,1)$}, \quad \mbox{and} \quad \frac{1}{4x}<\varphi(x)<\frac{1}{x}, \quad \mbox{for $x >1$.}
\end{equation}
In terms of $\varphi$, we may express $\sigma_f^2$ as
\[
 \sigma_f^2(t){=\sum_{k=1}^\infty \frac{b_k t}{(t+b_k)^2}} =\sum_{k=1}^\infty \varphi(t/b_k)\,, \quad \mbox{for any $t >0$.}
\]

Fix $n \ge 2$. If $t \le b_{n+1}$, we have, using \eqref{eq:estimates phi} and the growth condition~\eqref{eq:condition on bk}, that
\[
\sum_{k=n+1}^\infty  \varphi(t/b_k)<\sum_{k=n+1}^\infty \frac{t}{b_k}
<\frac{t}{b_{n+1}}\sum_{k=n+1}^\infty \frac{1}{2^{k{-}(n+1)}}=\frac{2t}{b_{n+1}}\cdot
\]
Analogously, if $t \ge b_{n-1}$,
\[
\sum_{k=1}^{n-1} \varphi(t/b_k)<\sum_{k=1}^{n-1} \frac{b_k}{t}<\frac{b_{n-1}}{t}\sum_{k=1}^{n-1} \frac{1}{2^{n-1{-}k}}<\frac{2 b_{n-1}}{t}\cdot
\]

%
%
%

We then have, for $n \ge 2$, that for any $t \in I_n$,
\[
\sum_{k=n+1}^\infty  \varphi(t/b_k)< 2 \,\sqrt{{b_n}/{b_{n+1}}} \qquad  \mbox{and}\qquad \sum_{k=1}^{n-1}  \varphi(t/b_k)<2\,\sqrt{{b_{n-1}}/{b_{n}}}.
\]

%

For $n\ge 2$, it follows then  that, if $t \in I_n$,
\begin{equation*}
\varphi(t/b_n)\le \sigma_f^2(t)\le \varphi(t/b_n)+4 \max\big\{\sqrt{{b_n}/{b_{n+1}}},\,\sqrt{{b_{n-1}}/{b_{n}}}\big\}.
\end{equation*}
The bounds in \eqref{eq:bounds sigma productos canonicos uniform} now follow by observing that
$\varphi(x)$ decreases whenever $x$ moves away from $1$.
%
\end{proof}

\subsection{Derivative power series \texorpdfstring{$\mathcal{D}_f$}{Df} and its family}
Let $f(z)=\sum_{n=0}^\infty a_n z^n$ be a power series in $\K$, with radius of convergence $R >0$,
and let $(X_t)_{t\in[0,R)}$ be  its associated family. We now consider the power
series $\mathcal{D}_f$ given by
\[
\mathcal{D}_f(z)=z f'(z)=\sum_{n=1}^\infty n a_n\, z^n, \quad \text{for $|z|<R$}.
\]
This power series $\mathcal{D}_f$ has radius of convergence $R$, but it is not in $\mathcal{K}$, as $\mathcal{D}_f(0)=0$.
In any case, we denote by $(W_t)_{t\in[0,R)}$ the associated family of random variables.

Observe that for  $t\in(0,R)$ and $n \ge 1$, we have that
\begin{equation}
\label{eq:probabilities of derivative}
\P(W_t=n)=\frac{n a_n t^n}{t f^{\prime}(t)}=\frac{1}{m_f(t)}\,\frac{n a_n t^n}{f(t)}=\frac{n}{m_f(t)}\,\P(X_t=n).
\end{equation}
Thus, for  $t\in(0,R)$, we have that
\[
m_{\mathcal{D}_f}(t)=\E(W_t)=\frac{1}{m_f(t)}\sum_{n=1}^\infty n^2\, \P(X_t=n)= \frac{\E(X_t^2)}{m_f(t)}=\frac{\E(X_t^2)}{\E(X_t)}\cdot
\]
If the power series $f$ has only two nonzero coefficients, $a_0$ and $a_N$, with $N \ge 1$, then
the random variables $W_t$ are constant, and ${\E(X_t^2)}/{\E(X_t)}=N$, for any $t \in (0,\infty)$ (in this special case,  $f$ is a polynomial and $R=+\infty$).
Otherwise (if $f$ has at least 3 nonzero coefficients), the quotient ${\E(X_t^2)}/{\E(X_t)}$ is
monotonically increasing in the interval $(0,R)$.

In general, we have that
\begin{equation}
\label{eq:moments of derivative}\E(W_t^p)=\frac{1}{m_f(t)} \, \E(X_t^{p + 1}), \quad \text{for any $p>0$ and any $t \in (0,R)$}.
\end{equation}

The quotient
\begin{equation}
\label{eq:quotient mDF by mf}
\frac{m_{\mathcal{D}_f(t)}}{m_f(t)}=\frac{\E(X_t^2)}{\E(X_t)^2}=\frac{\sigma_f^2(t)}{m_f^2(t)}+1,
\end{equation}
plays a relevant role in what follows.

\section{Growth of the moments of power series distributions}

Let $f$ be a power series in $\K$ with radius of convergence $R>0$ and with associated family $(X_t)_{t \in [0,R)}$.

\subsection{Growth and range of the mean \texorpdfstring{$m_f$}{mf}}
Since $X_t$ is not constant for any $t \in (0,R)$, we have that $\sigma_f^2(t)>0$, for any $t \in (0,R)$ and hence,
because of \eqref{eq:m and sigma in terms of f}, $m_f(t)$ is \emph{strictly increasing in $[0,R)$}, though, in general,
$\sigma_f(t)$ is not increasing.

We denote
\begin{equation}\label{eq:def de Mf}
M_f=\lim_{t \uparrow R} m_f(t).
\end{equation}
As recorded in the following {result}, 
 it is the case that $M_f=\infty$, except for some exceptional instances.

\begin{lemma}
[Lemma 2.2 of \cite{CFFM1}]
\label{lemma:caracterization of MF finito}
For $f(z)=\sum_{n=0}^\infty a_n z^n$ in $\mathcal{K}$ with radius of convergence $R>0$, we have $M_f<\infty$ in just the
following two cases\textup{:}
\begin{enumerate}
\item[\textup{(1)}] if $R<\infty$ and $\sum_{n=0}^\infty n a_n R^n <\infty$,
\item[\textup{(2)}] and if $R=\infty$ and $f$ is a polynomial.
\end{enumerate}

In the first case, we have $M_f=(\sum_{n=0}^\infty n a_n R^n) / (\sum_{n=0}^\infty a_n R^n)$.
For a polynomial $f\in \mathcal{K}$, we have $M_f=\mbox{\rm deg}(f)$.
\end{lemma}

As an incremental quotient, $m_f(t)$ and $\ln f(t)$ are related as in the following result.

\begin{lemma}
[Simi{\'c}, \cite{Simic}]
\label{lema:simic1}
For $\lambda >1$ and $t>0$, with $\lambda \, t<R$, we have
\[
m_f(t) \ln \lambda \le \ln \Big(\frac{f(\lambda t)}{f(t)}\Big)\le m_f(\lambda t) \ln \lambda.
\]
\end{lemma}

This Lemma \ref{lema:simic1} follows directly from the expression \eqref{eq:m and sigma in terms of f}
and the fact that $m_f$ is increasing.

\subsection{Relative growth of moments}

We are comparing now, \emph{in the case $M_f=+\infty$}, the growth of moments of different exponents
and also the growth of the factorial moments and the moments of the $X_t$; this is covered, respectively,
by Corollaries~\ref{cor:relative growth of moments} and~\ref{cor:growth of moment and factorial moment}.
The basic tool is the following lemma.

\begin{lemma}
\label{lema:comparison normalized moments}
For $1< \alpha < \beta$,
\begin{equation}
\label{eq:comparison normalized moments}
\frac{\E(X_t^\alpha)}{\E(X_t)^\alpha}\le \frac{\E(X_t^\beta)}{\E(X_t)^\beta} , \quad \text{for any $t\in (0,R)$.}
\end{equation}
\end{lemma}
\begin{proof} Recall that the $X_t$ are nonnegative variables.
Write $\alpha=u +(1-u)\beta$, with $u=(\beta-\alpha)/(\beta-1)\in (0,1)$.
H\"{o}lder's inequality gives that
\[
\E(X_t^\alpha)=\E(X_t^u \, X_t^{(1-u)\beta})\le \E(X_t)^u \E(X_t^\beta)^{1-u},
\]
and so that
\[
\frac{\E(X_t^\alpha)}{\E(X_t^\beta)}\le \Big(\frac{\E(X_t)}{\E(X_t^\beta)}\Big)^u.
\]
Jensen's inequality gives that $\E(X_t^\beta)\ge \E(X_t)^\beta$, and therefore,
\[
\frac{\E(X_t^\alpha)}{\E(X_t^\beta)}\le \E(X_t)^{(1{-}\beta)u}=\frac{\E(X_t)^{\alpha}}{\E(X_t)^\beta},
\]
as stated.
\end{proof}

From Lemma \ref{lema:comparison normalized moments}, we deduce the following two corollaries.

\begin{coro}\label{cor:relative growth of moments} Assume $M_f=+\infty$. If $1 < \alpha < \beta$, then
\[
\lim_{t\uparrow R} \frac{\E(X_t^\alpha)}{\E(X_t^\beta)}=0.
\]
\end{coro}
\begin{proof}
From Lemma \ref{lema:comparison normalized moments}, we have that
\[
\frac{\E(X_t^\alpha)}{\E(X_t^\beta)}\le \frac{\E(X_t)^\alpha}{\E(X_t)^\beta}=m_f(t)^{\alpha-\beta} , \quad \text{for any $t \in(0,R)$}.
\]
The statement follows since $\lim_{t \uparrow R}m_f(t)=+\infty$ and $\alpha -\beta<0$.
\end{proof}

\begin{coro}\label{cor:growth of moment and factorial moment}
Assume $M_f=+\infty$. For any integer $k \ge 1$, we have that
\[\lim_{t \uparrow R} \frac{\E(X_t^{\underline{k}})}{\E(X_t^k)}=1\, .
\]
\end{coro}
\begin{proof}From Corollary \ref{cor:relative growth of moments}, we deduce that
\[
\lim_{t \uparrow R} \frac{\E(X_t^{{j}})}{\E(X_t^k)}=0\,, \quad \text{for $0 \le j <k$.}
\]
The statement follows by expanding $\E(X_t^{\underline{k}})$ as $\E(X_t^k)$ plus a linear
combination of the moments $\E(X_t^j)$ with $0\le j <k$.
\end{proof}

\subsection{Growth and range of the variance \texorpdfstring{$\sigma_f^2$}{sigmaf\textasciicircum2}}\label{section:growth of variance}

Some of the results below concerning~$\sigma_f^2(t)$ will depend \emph{on the gaps among the indices} of the power series. We introduce next some convenient notation.

Let $(n_k)_{k\ge 1}$ be the increasing sequence of indices so that $a_{n_k}\neq 0$, for each $k \ge 1$, and $a_n=0$,
if $n \notin\{n_k: k \ge 1\}$. Thus the $(a_{n_k})$ are the nonzero Taylor coefficients of~$f$, and the random
variables $(X_t)$ take exactly the values $(n_k)$. Observe that $n_1=0$, and that for a polynomial $f$,  the
sequence $(n_k)_{k\ge 1}$ is finite.   Define $\text{gap}(f)$ and $\overline{\mbox{gap}}(f)$ as
\begin{equation}
\label{eq:def de gaps}
\mbox{gap}(f)=\sup_{k \ge 1} \,(n_{k + 1}-{n_k})\quad \text{and} \quad
\overline{\mbox{gap}}(f)=\limsup_{k \to \infty} \,(n_{k + 1}-{n_k}).
\end{equation}

It is always the case that $\mbox{gap}(f)\ge \overline{\mbox{gap}}(f) \ge 1$, for any $f \in \K$ which is not a polynomial.
For polynomials, we still have $\mbox{gap}(f)\ge 1$, and we can define $\overline{\mbox{gap}}(f)=0$.

We divide the discussion on variance growth depending on whether $R$ is infinite or finite.
\subsubsection{Entire functions, \texorpdfstring{$R=\infty$}{R=infinity}}\label{subsec:entire}

\mbox{}
\medskip

A. \emph{Lower bounds for $\sup_{t>0} \sigma_f^2(t)$ and for $\limsup_{t \to \infty} \sigma_f^2(t)$.}

\smallskip

The (universal) lower bounds on $\sigma_f$ that we are about to discuss originate with Hayman's results in~\cite{Hayman2} quantifying
Hadamard's three lines theorem for entire functions (not necessarily in $\K$). See also \cite{Abi2} and \cite{Kjellberg}.

\begin{theor}
[{Bo\u{\i}chuk}--Gol'dberg, \cite{BoichukGoldberg}]\label{teor:bound for sigma}
If\/ $f \in \K$ is entire, then
\[
\sup_{t>0} \sigma_f^2(t) \ge \frac{1}{4} \, \mbox{\rm gap}(f)^2\quad \Big(\ge \frac{1}{4}\Big).
\]
If, moreover, $f$ is transcendental, then
\[
{\limsup_{t\to\infty} \sigma_f^2(t) \ge \frac{1}{4} \, \overline{\textup{gap}}(f)^2\quad \Big(\ge \frac{1}{4}\Big)}.
\]
\end{theor}

Recall that an entire function $f$ is termed \emph{transcendental} if it is not a polynomial.

\smallskip

This result is Theorem 2 of \cite{BoichukGoldberg}. See also Theorem 1 {in} 
 \cite{Abi} and Lemma~2.5 in~\cite{OstrovskiiUreyen}.
The probabilistic argument below is simpler than the original proof; it uses the discreteness of the random
variables $X_t$.

\begin{proof}[Proof of Theorem {\textup{\ref{teor:bound for sigma}}}]
Let $(n_k)^N_{k=1}$ be the indices of the nonzero coefficients of $f$, with $N\le +\infty$.

Fix $k <N$ and take $t^\star>0$ so that $m_f(t^\star)=(n_{k + 1}+n_k)/2$, i.e., the midpoint of
the interval $[n_k, n_{k + 1}]$. Such $t^\star$ exists because $m_f(t)$ is a continuous (and increasing)
function, $m_f(0)=0$, and $M_f=\infty$ or $M_f=\text{degree}(f)$ if $f$ is a polynomial (recall
Lemma~\ref{lemma:caracterization of MF finito}).

As $X_{t^\star}$ takes the values $n_1,n_2,\dots$, clearly $|X_{t^\star}-m_f(t^\star)|\ge \frac{1}{2} (n_{k + 1}-n_k)$
with probability 1. This gives
\[
\sigma^2_f(t^\star) =\E\big((X_{t^\star}-m_f(t^\star))^2\big)\ge \frac{1}{4} (n_{k + 1}-n_k)^2.
\]
The statements now follows by taking sup and limsup in the inequality above and appealing
to the definitions of ${\textup{gap}}$ and $\overline{\textup{gap}}$.
\end{proof}

In fact, the very same argument shows, for the centered moments,  that if $f \in \K$ is entire, then 
\[
\sup_{t>0} \E (|X_t-m_f(t)| )^p \ge \frac{1}{2^p} \, \mbox{\rm gap}(f)^p,\quad \mbox{for any $p>0$},
\]
and,  moreover, that if $f$ is transcendental, then
\[
\limsup_{t\to\infty} \E (|X_t-m_f(t)| )^p \ge \frac{1}{2^p}\, \overline{\textup{gap}}(f)^p, \quad \mbox{for any $p>0$}.
\]

As for the \emph{sharpness of the sup part} of Theorem \ref{teor:bound for sigma}, consider the case $f(z)=a+bz$, with $a,b>0$,
for which $X_t$ is a Bernoulli variable with success probability $bt/(a+bt)$. In this case,
one has $\sigma^2_f(t)=abt/(a+bt)^2$, which takes its maximum value of $1/4$ at $t=a/b$.

In fact, the converse is also true.
\begin{theor}[Abi-Khuzzam, \cite{Abi}]\label{teor:abi}
For $f\in \K$ entire, $\sup_{t >0} \sigma_f^2(t)=1/4$ if and only if $f(z)=a+bz$, with $a,b >0$.
\end{theor}

This result is Theorem 3 of \cite{Abi}. See also Lemma 2.5 in~\cite{OstrovskiiUreyen}.
The probabilistic argument below is again simpler than the original proof.

\begin{proof}[Proof of Theorem {\textup{\ref{teor:abi}}}] The `if' part has been discussed above.

Assume that $f\in \K$ is entire and that $\sup_{t >0} \sigma_f^2(t)=1/4$. Theorem \ref{teor:bound for sigma} gives
that $\mbox{gap}(f)=1$. Let $a_n$ and $a_{n + 1}$ be any two nonzero consecutive coefficients of $f$, and let~$t^\star$ be
such that $m_f(t^\star)=n+1/2$ (observe that, in any case, $M_f \ge n + 1$). We have that $|X_{t^\star}-m_f(t^\star)|\ge 1/2$.
By hypothesis, $\E((X_{t^\star}-m_f(t^\star))^2)=\sigma_f^2(t^\star)\le 1/4$, and thus $|X_{t^\star}-m_f(t^\star)|= 1/2$ with probability 1,
which means that $X_{t^\star}$ only takes the values $n$ and $n + 1$, and thus $f(z)=a_n z^n+a_{n + 1} z^{n + 1}$.
Since $f(0)>0$, because $f \in \K$, it must be the case that $n=0$ and $f(z)=a+bz$, with $a,b>0$.
\end{proof}

For the \emph{sharpness of the limsup part} of Theorem \ref{teor:bound for sigma}, consider a canonical product $h$ given by the infinite product
\begin{equation}\label{eq:Hayman example}
h(z)=\prod_{n=1}^\infty \Big(1+\frac{z}{b_n}\Big),
\end{equation}
where $(b_n)_{n\ge 1}$ is a sequence of positive numbers increasing to $+\infty$ so that \eqref{eq:condition on bk} holds and, in fact, such that $\lim_{n \to \infty} b_{n+1}/b_n=+\infty$; in particular, $\sum_{k\ge 1} 1/b_k<+\infty$ holds. Obviously,  $\mbox{gap}(h)=1$ and $\overline{\mbox{gap}}(h)=1$. It {follows directly from the estimates \eqref{eq:bounds sigma productos canonicos uniform} that $\limsup_{t \to \infty}\sigma^2_h(t)=1/4$; this is Hayman's example in Theorem 4 of~\cite{Hayman2}.} If $h$ is multiplied by a polynomial $p$ in such a way
that $f=ph \in \K$, then it is still the case that $\limsup_{t \to \infty}\sigma^2_f(t)=1/4$.

As it turns out, Abi-Khuzzam has characterized {(see Theorem 2 in \cite{Abi2} and its proof)} the entire
functions $f\in \K$ with $\limsup_{t\to \infty} \sigma_f^2(t)=1/4$ as precisely those entire
functions $f\in \K$ which factorize as
\[
f(z)=p(z) \prod_{n=1}^\infty \Big(1+\frac{z}{b_n}\Big),
\]
where the $b_n$ are as in Hayman's example, and where $p$ is a polynomial.

\medskip

B. \emph{Limit of\/ $\sigma^2_f(t)$ as $t \to \infty$.}

\smallskip

Regarding  the existence and possible limits of $\sigma_f^2(t)$ as $t \to \infty$, the following holds.

\smallskip
{B.1.} \emph{Polynomials.} Polynomials $f \in \K$ are characterized, among the entire functions in~$\K$, by
\begin{equation}\label{eq:limit sigma for polys}\lim_{t \to \infty} \sigma_f^2(t)=0.\end{equation}
A direct calculation with the formulas \eqref{eq:m and sigma in terms of f} shows that for polynomials \eqref{eq:limit sigma for polys} holds; in fact $\sigma_f^2(t)=O(1/t)$, as $t \to \infty$. The converse follows, for instance, from Theorem \ref{teor:bound for sigma}.

\smallskip

{B.2.} \emph{Transcendental {functions}.}
As shown by Hilberdink in \cite{Hilberdink},
for a transcendental entire function \emph{it is never the case that $\lim_{t \to \infty} \sigma_f(t)$ exists and it is finite}.

{However}, 
 there are entire functions  $f \in \K$  for which $\limsup_{t \to \infty} \sigma_f^2(t)<+\infty$ and $\liminf_{t\to \infty} \sigma_f^2(t)>0$. To see this, just consider a canonical product $f$ as in \eqref{eq:formula of canonical product}, with $b_n=2^n$, and apply \eqref{eq:bounds sigma productos canonicos uniform}.

It is {also} possible to have  $\lim_{t \to \infty} \sigma_f(t)=\infty$, as shown, for instance, by the
exponential function $f(z)=e^z$, where $\sigma_f^2(t)=t$, for $t \ge 0$.

We emphasize that for $f \in \K$ entire, if $\lim_{t \to \infty}\sigma_f^2(t)$ exists, then that limit is $0$ (just for polynomials) or $+\infty$.

\medskip

C. \emph{Boundedness of\/ $\sigma_f^2(t)$.}

\smallskip

We discuss now when, if ever,
\begin{equation}
\label{eq:sup sigma finito}
 \sup_{t>0} \sigma_f(t)<+\infty
\end{equation}
does 
 hold for an entire function $f \in \K$.

 For polynomials $f\in \K$, \eqref{eq:sup sigma finito} holds since, in fact, in this case, $\lim_{t \to \infty}\sigma_f(t)=0$.

In general, \emph{if\/ \eqref{eq:sup sigma finito} holds, then the order $\rho(f)$ of $f$ must be zero.}
(See Section~\ref{section:entire functions} for details about the order of an entire function $f$ in $\mathcal{K}$.)

To see this, observe that, since $tm_f^\prime(t)=\sigma_f^2(t)$, integrating, we deduce that $m_f(t)=O(\ln t)$, as $t \to \infty$.
A further integration, using \eqref{eq:m and sigma in terms of f},  shows that
\begin{equation}
\label{eq:asymp of ln f}
\ln f(t)=O((\ln t)^2)\,, \quad \text{as $t \to \infty$.}
\end{equation}
This gives that $\rho(f)=0$, see \eqref{eq:formula order}.

Alternatively, Proposition \ref{prop:orden y sigma cuad partod por m} below shows that
\[
\rho(f)\le \big(\sup_{t>0} \sigma^2_f(t)\big)\frac{1}{M_f}\cdot
\]
If $f$ is a polynomial, then $\rho(f)=0$; and if $f$ is transcendental, $M_f=\infty$, and thus $\rho(f)=0$.

\smallskip

However,
  $\rho(f)=0$, or even the stronger condition \eqref{eq:asymp of ln f}, {are} 
  not enough to ensure the boundedness of $\sigma_f^2(t)$.

Consider the canonical product
\[
g(z)=\prod_{k=1}^\infty \Big(1+\frac{z}{b_k}\Big)^k,
\]
where $(b_k)_{k\ge 1}$ is a sequence of positive numbers increasing to $\infty$ satisfying \eqref{eq:condition on bk} and such that $\sum_{k=1}^\infty k/b_k<+\infty$.
{For each $k\ge 1$,  $-b_k$ is a zero with multiplicity $k$ of the entire function $g \in \K$}.

{Additionally,} we assume also that for some constant $H>0$, the $b_k$ satisfy
\[
\sum_{k<n}k b_k <H b_n \quad \mbox{and} \quad \sum_{k>n} \frac{k}{b_k}\le \frac{H}{b_n} \, \quad \mbox{for each $n \ge 2$.}
\]

In this case,
\[
\sigma_g^2(t)=\sum_{k=1}^\infty k \,\varphi\Big(\frac{t}{b_k}\Big)\, , \quad \mbox{for any $t >0$,}
\]
{where $\varphi(x)={x}/{(1+x)^2}$}.

Let $C$ denote a generic positive constant. With the notations of Section \ref{section:mean variance canonical products} and estimating as in there, we obtain
\[
n \,\varphi\Big(\frac{t}{b_n}\Big)\le \sigma_g^2(t) \le n +C\, n \max\big\{\sqrt{{b_{n-1}}/{b_n}}, \sqrt{{b_{n}}/{b_{n+1}}}\,\big\}\, , \quad \mbox{for $t \in I_n$ and $n \ge 2$,}
\]
where $I_n=[\sqrt{b_{n-1}b_n}, \sqrt{b_n b_{n+1}}]$.
It follows, in particular, since $\sigma_g^2(b_n)\ge n/4$,  that $\limsup_{t \to \infty} \sigma_g^2(t)=+\infty$, and also that,
\[
\sup_{t \in I_n} \sigma_g^2(t) \le  C n.
\]
From \eqref{eq:variance less than mean of canonical products}, we see that
\begin{equation}\label{eq:bound for mg special can prods} m_g(t)\le C n +\sum_{k\le n}k\le C n^2\, \quad \mbox{for any $t\in I_n$ and $n \ge 2$.}\end{equation}

Consider now the specific sequence $b_k=e^{e^k}$, $ k \ge 1$, which satisfies the requirements above.
In this case, \eqref{eq:bound for mg special can prods} translates into
\[
m_g(t)\le C (\ln\ln t)^2\,, \quad \mbox{for $t\ge 2$}.
\]
Integrating, this gives, for this  example, that
\[
\ln g(t)=O(\ln t (\ln\ln t)^2)\, \quad \mbox{as $t \to \infty$,}
\]
which implies $\rho(g)=0$. Notice that for any function $\Phi(t)$ slowly increasing to $\infty$, the sequence $(b_k)_{k\ge 1}$ can be chosen so that $\ln g(t)=O(\Phi(t) \ln t)$ just by making $b_k$ increase fast enough. This is the best that can be expected, {because} if  an entire function $h$ in $\K$ grows as $\ln h(t)=O(\ln t)$ as $t \to \infty$, then $h$ is a polynomial.

\subsubsection{Finite radius:  \texorpdfstring{$R<\infty$}{R<infinity}}

\mbox{}
\medskip

We now turn to functions $f \in \K$ with finite radius $R$ of convergence. We have the following results
on the behaviour of $\sigma^2_f(t)$. We divide the discussion according to whether~$M_f$  is finite or not.

\vskip3pt

$\bullet$ \emph{Case $M_f=+\infty$.}

\vskip2pt

In this case,
\begin{equation}\label{eq:sup sigma2 infinity} \sup_{t \in (0,R)} \sigma^2_f(t) =+\infty\,,
\end{equation}
and also, of course, $\limsup_{t \uparrow R} \sigma^2_f(t)=+\infty$.
To verify \eqref{eq:sup sigma2 infinity}, assume that $\sup_{t \in (0,R)} \sigma^2_f(t)=S<+\infty$. Thus $tm_f^\prime(t)\le S$, for $t \in [0,R)$.
Integrating between $R/2$ and $t\in (R/2,R)$, we would have that
\[
m_f(t)\le m_f(R/2) + S \ln \Big(\frac{2t}{R}\Big) , \quad \text{for $t \in (R/2, R)$,}
\]
which implies, by letting $t \uparrow R<+\infty$, that
$M_f\le m_f(R/2)+S \ln 2<+\infty$.

\vskip3pt

{$\bullet$ \emph{Case $M_f<+\infty$.}}

\vskip2pt

 If $M_f<+\infty$, then $\sum_{n=0}^\infty n a_n R^n <+\infty$,
see Lemma \ref{lemma:caracterization of MF finito}, and in fact
\[
\Sigma:=\lim_{t \uparrow R} \sigma_f^2(t)=\frac{\sum_{n=0}^\infty n^2 a_n R^n}{\sum_{n=0}^\infty a_n R^n}
-\bigg(\frac{\sum_{n=0}^\infty n a_n R^n}{\sum_{n=0}^\infty a_n R^n}\bigg)^2.
\]

It is always the case that $\lim_{t \uparrow R} \sigma_f^2(t) >0$, since $\Sigma$ is the variance of the
random variable~$Z$ that takes, for each integer $n \ge 0$, the value $n$ with probability
$a_n R^{n}/(\sum_{k=0}^\infty a_k R^k)$, and $Z$ is a nonconstant variable since $f$ is in $\K$.

But there is no absolute positive lower bound for $\lim_{t \uparrow R} \sigma_f^2(t)$.
For $\varepsilon >0$, the power series $f(z)=1+\varepsilon \sum_{n=1}^\infty z^n/n^4$ is in $\K$ and has radius of convergence $R=1$.
We have
\[
M_f=\frac{\varepsilon \zeta(3)}{1+\varepsilon \zeta(4)} \quad\text{and}\quad
\lim_{t \uparrow 1} \sigma_f^2(t)=\frac{\varepsilon \zeta(2)}{1+\varepsilon \zeta(4)}
-\Big(\frac{\varepsilon \zeta(3)}{1+\varepsilon \zeta(4)}\Big)^2,
\]
which tends to 0 as $\varepsilon \downarrow 0$.

\smallskip

The example  $f(z)=\sum_{n=0}^\infty z^n/(1 + n)^3$ shows that  $\lim_{t \uparrow R} \sigma_f^2(t)=\infty$ may happen.

\subsection{Growth of the quotient \texorpdfstring{$\sigma_f/m_f$}{sigmaf/mf} and gaps}
\label{section:growth of quotient}

Let $f \in \K$, not a polynomial,
have radius of convergence $R>0$. As in Section~\ref{section:growth of variance},
we denote by $(n_k)_{k=1}^\infty$ the increasing sequence of indices of the nonzero
coefficients of $f$. We define $\overline{G}(f)$ by
\begin{equation}
\label{eq:def de G barra}
\overline{G}(f)=\limsup_{k\to \infty} \,\frac{n_{k + 1}}{n_k}\cdot
\end{equation}
Clearly, $\overline{G}(f)\ge 1$. If for a given $f$ we had $\liminf_{k\to \infty} n_{k + 1}/{n_k}>1$,
then $f$ would be the sum of a polynomial and a power series with Hadamard gaps; but notice that $\overline{G}(f)$
calls for a \lq$\limsup$\rq. Observe that the definition of $\overline{G}(f)$ involves the quotients $n_{k + 1}/{n_k}$,
and not the differences $n_{k + 1}-n_{k}$ as it is the case in $\mbox{gap}(f)$ and $\overline{\mbox{gap}}(f)$.

\begin{theor}
\label{toer:limsup del quotient}
Assume that $f\in \K$, with radius of convergence $R>0$, is not a polynomial and that $M_f=+\infty$. Then
\[
\limsup_{t\uparrow R} \frac{\sigma_f(t)}{m_f(t)}\ge \frac{\overline{G}(f)- 1}{\overline{G}(f)+ 1}\cdot
\]
\end{theor}

The proof below mimics our proof of Theorem \ref{teor:bound for sigma}.

\begin{proof}
Since $M_f=+\infty$, we have that $m_f(t)$ is a homeomorphism from $[0,R)$ onto $[0,+\infty)$,
and thus for any integer $k$ there is $t_k \in (0,R)$ so that
\[
m_f(t_k)=\frac{n_k + n_{k+1}}{2}\cdot
\]
For the random variable $X_{t_k}$, we have that
$|X_{t_k}-m_f(t_k)|\ge (n_{k+1}-n_k)/2$ with probability~$1$,
and thus
\[
\sigma^2_f(t_k)=\E\big((X_{t_k}-m_f(t_k))^2\big)\ge \frac{1}{4} \,(n_{k+1}-n_k)^2,
\]
and also
\[
 \frac{\sigma^2_f(t_k)}{m^2_f(t_k)}\ge \frac{(n_{k+1}-n_k)^2}{(n_{k+1} + n_k)^2}\cdot
\]
The result follows.
\end{proof}

Similarly, for the general centered moments of the family of functions $f\in \K$ as in the statement of Theorem \ref{toer:limsup del quotient}, we have that
\[
\limsup_{t\uparrow R} \frac{\E\left(|X_t-m_f(t)|^p\right)}{m_f(t)^p}\ge \Big(\frac{\overline{G}(f)- 1}{\overline{G}(f)+ 1}\Big)^p, \quad \mbox{for any $p>0$}.
\]

\subsection{Zero-free region and \texorpdfstring{$\sigma_f$}{sigmaf}}

A function $f$ in $\mathcal{K}$ does not vanish on the interval $[0,R)$ and, in fact, its zeros
must lie away from that segment.

The following result shows how the variance function of $f$ determines a specific zero free region
for $f\in\mathcal{K}$ containing the interval $[0,R)$.

\begin{prop}\label{prop:ceros de funcion f de K}
Let $f\in \mathcal{K}$ have radius of convergence $R>0$.
If for some $t\in [0,R)$ and some $\theta \in [-\pi, \pi]$, we have $f(te^{\imath \theta})=0$, then
\[
|\theta|\cdot\sigma_f(t)\ge \frac{\pi}{2}\cdot
\]
Thus, $f\in\mathcal{K}$ does not vanish in the region
\[
\Omega_f=\Big\{z=t e^{\imath \theta}: t \in[0,R) \mbox{ and } |\theta|< \frac{\pi}{2\sigma_f(t)}  \Big\}.
\]
\end{prop}
For the proof, we may use the following lemma.
\begin{lemma}
[Sakovi{\v{c}}, \cite{Sakovic}]
\label{lema:ceros caracteristica}
Let $Y$ be a random variable and let $\theta \in \R$. If\/ $\E(e^{\imath \theta Y})=0$, then
\[
\theta^2 \,\V(Y)\ge \dfrac{\pi^2}{4}\cdot
\]
\end{lemma}

The bound on Lemma \ref{lema:ceros caracteristica} is sharp. Simply consider the random variable $Z$
which takes values $\pm1$ with probability $1/2$; then $\V(Z)=1$ and $\E(e^{\imath \theta Z})=\cos \theta$,
which vanishes at~$\pi/2$. (Actually, equality in Lemma \ref{lema:ceros caracteristica} only happens
for this simple symmetric random variable $Z$.)

The result of Sakovi{\v{c}} appeared in \cite{Sakovic}. As presented in the more
accesible reference~\cite{Rossberg}, Lemma \ref{lema:ceros caracteristica} follows most ingeniously as follows.

\begin{proof}[Proof of Lemma \textup{\ref{lema:ceros caracteristica}} following Rossberg \cite{Rossberg}]
Consider the function
\[
\varphi(t)=t^2-1+\frac{4}{\pi} \cos \frac{\pi t}{2},
\]
which happens to be positive for all $t \in \R$, except for $t=\pm 1$, where $\varphi(t)=0$.
(A~misprinted sign in the definition of $\varphi$ in~\cite{Rossberg} has been corrected.)
Assume that $\E(e^{\imath \theta Y})=0$. Consider $W=Y-\E(Y)$, so that $\E(e^{\imath \theta W})=0$,
and thus $\Re \E(e^{\imath \theta W})=\E(\cos(\theta W))=0$.
This yields
\[
0\le \E\Big(\varphi\Big( \frac{2\theta}{\pi} \,W\Big)\Big)=\theta^2 \,\frac{4}{\pi^2} \,\E(W^2)-1
=\theta^2 \,\frac{4}{\pi^2} \,\V(Y)-1.\qedhere
\]
\end{proof}

A more direct proof of Lemma \ref{lema:ceros caracteristica}, but with a weaker constant, appears, for instance,
in B\'aez-Duarte~\cite{BaezDuarteOtro}, see Proposition 7.8 (see also Lemma~2.3 in~\cite{CFFM1}).

\begin{proof}
[Proof of Proposition \textup{\ref{prop:ceros de funcion f de K}}]
It follows from Lemma \ref{lema:ceros caracteristica} and from observing that if $f(te^{\imath \theta})=0$,
then $\E(e^{\imath \theta X_t})=f(t e^{\imath \theta})/f(t)=0$.
\end{proof}

Alternatively, to verify Proposition \ref{prop:ceros de funcion f de K}, we may use
Lemma 1 of \cite{Abi}, which gives that for $f \in \K$,
\[
f(t)^2-|f(te^{\imath \theta})|^2\le 4 \sin^2(\theta/2) \,f(t)^2 \,\sigma^2_f(t),
\quad \text{ for $t\in (0,R)$ and $\theta \in [-\pi, \pi]$}.
\]
If $f(te^{\imath \theta})=0$, then $1\le 2 |\sin(\theta/2)| \,\sigma_f(t)$, and
thus $|\theta|\,\sigma_f(t)\ge 1$. This gives a weaker result with the constant $\pi/2$
replaced by 1. Lemma 1 of \cite{Abi} is stated for entire functions, but it is valid
for general $f \in \K$.

\begin{remark}
\label{remark:BoichukGoldberg from Sakovic}
(\emph{On {Bo\u{\i}chuk}--Gol'berg's Theorem \textup{\ref{teor:bound for sigma}}}).
From Sakovi{\v{c}}'s Lemma \ref{lema:ceros caracteristica}, we may deduce
{Bo\u{\i}chuk}--Gol'berg's Theorem \ref{teor:bound for sigma} in the weaker form that
\[
\sup_{t >0}\sigma_f^2 (t)\ge 1/4,
\]
for any entire function $f \in \K$. To see this, we may assume that $\sup_{t >0}\sigma_f^2 (t)<+\infty$.
As discussed in Subsection~\ref{subsec:entire}, this yields that the entire function $f$ is of order $0$,
and Hadamard's factorization theorem gives that~$f$ is an infinite canonical product or a polynomial
(nonconstant and not a monomial, since $f\in \K$). (This is the starting point of Hayman's proof of Theorem 3
of \cite{Hayman2}.) In any case, $f$ vanishes at some $z_0\neq 0$. Write $z_0=r_0\, e^{\imath \theta_0}$,
with $r_0>0$ and $|\theta_0|\le \pi$. Lemma \ref{lema:ceros caracteristica} gives that
\[
\pi \sigma_f(r_0)\ge \frac{\pi}{2},
\]
and thus, that $\sigma_f(r_0)\ge 1/2$, and, in particular, that $\sup_{t >0}\sigma_f^2 (t)\ge 1/4$.

This same reasoning also gives that if $f$ is not a polynomial, then
\[
\limsup_{t \to \infty}\sigma_f^2 (t)\ge 1/4.
\]
\end{remark}

\section{Clans}\label{section:clans}

In previous sections, we have compared the growth of $\E(X_t^\beta)$ with $\E(X_t)^\beta$, as~$t\uparrow R$,
for the basic examples of Khinchin families. Motivated by Hayman \cite{Hayman} (see Remark~\ref{remark:momentos de Hayman}),
we introduce next a particular kind of Khinchin families which we shall call \emph{clans},
for which $\E(X_t^2)\sim \E(X_t)^2$ as $t\uparrow R$. Concretely,

\begin{defin}
\label{def:clan}
Let $f$ be in $\mathcal{K}$ have radius of convergence $R\le \infty$.
We say that $f$ is a \emph{clan} \textup{(}and also that the associated family $(X_t)_{t \in [0,R)}$
is a clan\textup{)} if
\begin{equation}
\label{eq:lim para clan}
 \lim_{t \uparrow R}\frac{\sigma_f(t)}{m_f(t)}=0.
\end{equation}
\end{defin}

For a clan, the normalized variables $Y_t=X_t/\E(X_t)$, for $t \in (0,R)$, converge in probability
to the constant 1 as $t \uparrow R$, since its variance $\V(Y_t)={\sigma_f^2(t)}/{m_f(t)^2}$ converges
to~0 as $t \uparrow R$.

This clan condition is equivalent to
\[
\lim_{t \uparrow R} \frac{\E(X_t^2)}{\E(X_t)^2}=1.
\]
In terms of just the mean $m_f$, being a clan, see \eqref{eq:quotient mDF by mf},  is equivalent to
\[
\lim_{t \uparrow R} \frac{t m_f^\prime(t)}{m_f(t)^2}=0,
\]
while in terms of the derivative power series $\mathcal{D}_f$, the condition for being clan becomes
\[
\lim_{t \uparrow R} \frac{m_{\mathcal{D}_f}(t)}{m_f(t)}=1.
\]
Alternatively, if we define
\begin{equation}
\label{eq:def de Lf}
L_f(t):=\frac{f(t)f^{\prime\prime}(t)}{f'(t)^2}, \quad \text{for $t \in (0,R)$},
\end{equation}
and since
\begin{equation}
\label{eq:clan in terms of f previo}\frac{\E(X_t^2)}{\E(X_t)^2}
=\frac{1}{m_f(t)}+L_f(t),
\end{equation}
we have that $f$ is a clan if and only if
\begin{equation}
\label{eq:clan in terms of f}
\lim_{t \uparrow R}L_f(t)=1-{1}/{M_f}\cdot
\end{equation}

\subsection{Some examples}

\subsubsection{Examples of clans}\label{seccion:examples of clans}

The Poisson family associated to $f(z)=e^z$ is a clan, since in this case $m_f(t)=t$ and $\sigma_f^2(t)=t$,
and the radius of convergence is $R=\infty$.

The Bernoulli and binomial families are also clans. In fact, all polynomials $f$ are clans,
since for them $\sigma_f(t)\to 0$, while $m_f(t) \to \textup{degree}(f)$, as $t \to \infty$.
Further, we have the following.

\begin{lemma}
\label{lem:polynomial clans}
Let $f\in \mathcal{K}$ be a clan. Then $M_f<+\infty$ if and only if $f$ is a polynomial.
\end{lemma}

\begin{proof}
Polynomials in $\K$ are clans and have $M_f$ finite; in fact, $M_f$ coincides with its degree.

To show the converse, let $f(z)=\sum_{n=0}^\infty a_n z^n$ be in $\mathcal{K}$, not a polynomial, and with $M_f<\infty$.
Then, because of Lemma \ref{lemma:caracterization of MF finito}, we have that $R<\infty$ and
$\sum_{n=1}^\infty n a_n R^n< +\infty$.

Then since
\[
\frac{\E(X_t^2)}{\E(X_t)^2}=\frac{(\sum_{n=0}^\infty n^2 a_n t^n) (\sum_{n=0}^\infty a_n t^n)}{(\sum_{n=0}^\infty n a_n t^n)^2},
\quad \text{for $t \in (0,R)$},
\]
and since $f$ is a clan,
taking limit as $t\uparrow R$, we obtain that
\[
\frac{(\sum_{n=0}^\infty n^2 a_n R^n) (\sum_{n=0}^\infty a_n R^n)}{(\sum_{n=0}^\infty n a_n R^n)^2}=1.
\]
If we now set $b_n=a_nR^n/(\sum_{j=0}^\infty a_j R^j)$, for each $n \ge 0$, which satisfy $\sum_{n=0}^\infty b_n=1$,
then the identity above becomes
\[
\sum_{n=0}^\infty n^2 b_n=\Big(\sum_{n=0}^\infty n b_n\Big)^2 .
\]
This means that $b_n=0$, for each $n \ge 0$ except for one value of $n$, which would imply
that $f$ is a monomial: a contradiction.
\end{proof}

The partition function $P(z)=\prod_{k=1}^\infty 1/(1-z^k)$ and the Bell function $B(z)=e^{e^z-1}$
are also clans. This follows immediately from \eqref{eq:mean and variance partition}
and \eqref{eq:mean and variance Bell}.

{In  Section~\ref{section:Strongly Gaussian Khinchin families}, we exhibit an ample class of functions, that includes the generating
function of the partitions and its variants, which are clans}.

\subsubsection{Clans and $L_f$}

The characterization of clan given in \eqref{eq:clan in terms of f} using the function $L_f$ gives
immediately the following.

\begin{lemma}
\label{lem:condition for clan if Mfinfinite}
Let $f$ be a power series in $\K$ with radius of convergence $R\le +\infty$ and such that $M_f=+\infty$.
Then $f$ is a clan if and only if
\begin{equation}
\label{eq:condition for clan if Mfinfinite R}
\lim_{t \uparrow R} L_f(t)=1.
\end{equation}
In particular, if $f$ is an \emph{entire transcendental function in $\K$}, then $f$ is a clan if and only
if \eqref{eq:condition for clan if Mfinfinite R} holds.
\end{lemma}

Recall, from Lemma~\ref{lemma:caracterization of MF finito}, that entire transcendental
functions in $\K$ have $M_f=\infty$. For any transcendental entire function, it is always the case that
\begin{equation}\label{eq:liminf Lf ge 1}
\liminf_{t \to \infty} L_f(t)\ge1.
\end{equation}
This follows from Lemma \ref{lemma:caracterization of MF finito} and the general identity \eqref{eq:clan in terms of f previo}.
This has been pointed out by Simi{\'c} in \cite{Simic0}, p.~682. {But, in fact, we have the following}.

\begin{lemma}\label{lema:previo weak clan} For any transcendental entire function $f$ in $\K$,
\begin{equation}\label{eq:liming L=1, transcendental}\liminf_{t \to \infty} L_f(t)=1.
\end{equation}
\end{lemma}

\begin{proof}
Since $f$ is trascendental, we have that $M_f=+\infty$.
Let $c$ be such that $m_f(c)=1$. From the identity
\[L_f(s)=\Big(s\Big(1-\frac{1}{m_f(s)}\Big)\Big)^\prime\, , \quad \mbox{for $s >0$,}
\]
we deduce, since $M_f=+\infty$,  that
\[
\int_{c}^{t} L_f(s) \,ds=t \Big(1-\frac{1}{m_f(t)}\Big)\, , \quad \mbox{for $t >c$,}
\]
and so, that
\[\lim_{t \to \infty}\frac{1}{t}\int_{c}^{t} L_f(s) ds=1.
\]
This implies, given \eqref{eq:liminf Lf ge 1}, that~\eqref{eq:liming L=1, transcendental} holds.
\end{proof}

\subsubsection{Hayman class, strong gaussianity and clans}
\label{section:Strongly Gaussian Khinchin families}

Let $f$ be in $\K$ with radius of convergence $R>0$, and let $(X_t)_{t \in [0,R)}$ be the associated family.
Write $\breve{X}_t$ for the normalized random variable
\[
\breve{X}_t=\frac{X_t-m_f(t)}{\sigma_f(t)}\, \quad \text{for $t \in (0,R)$}.
\]
The characteristic function of $\breve{X}_t$ is
\[
\E(e^{\imath \theta \breve{X}_t})=\E(e^{\imath \theta X_t/\sigma_f(t)}) \, e^{-\imath \theta m_f(t)/\sigma_f(t)},
\quad \text{for $t\in (0,R)$ and $\theta \in \R$}.
\]
As introduced by B\'aez-Duarte in \cite{BaezDuarte}, a power series $f \in \K$ and its family $(X_t)_{t \in [0,R)}$
are termed \emph{strongly Gaussian} if
\[
\lim_{t \uparrow R} \sigma_f(t)=+\infty,
\quad \text{ and }\quad
\lim_{t \uparrow R} \int\nolimits_{|\theta|<\pi \sigma_f(t)}\big|\E(e^{\imath \theta \breve{X}_t})-e^{-\theta^2/2}\big| \, d\theta=0.
\]
Every strongly Gaussian function is Gaussian, that is, its normalized Khinchin family, $(\breve{X}_t)_{t\in[0,R)}$
converges in distribution, as $t\uparrow R$, to the standard normal variable or, equivalently,
\[
\lim_{t \uparrow R} \E(e^{\imath \theta \breve{X}_t})=e^{-\theta^2/2},
\quad \text{for each $\theta \in \R$}.
\]
See \cite{CFFM1} for definitions and proofs. The main interest of strong gaussianity
is that if a power series $f(z)=\sum_{n\ge 0}a_nz^n$ in $\K$ is strongly Gaussian,
then
\begin{equation}
\label{eq: asintotico coeficientes}
a_n\sim \frac{f(t_n)}{\sqrt{2\pi}\,t_n^n\,\sigma_f(t_n)},\quad \text{as $n\to\infty$},
\end{equation}
where $t_n$ is the unique value such that $m_f(t_n)=n$.

All \emph{\textup{(}Hayman\textup{)} admissible functions}, in the terminology of~\cite{Hayman},
see Definition in pages 68--69, or functions in the \emph{Hayman class}, are strongly Gaussian.
See, for instance, Theorem 3.8 in \cite{CFFM1}.

The exponential $f(z)=e^z$ is strongly Gaussian. In \cite{CFFM1} and \cite{CFFM2}, some criteria
are given to check when a power series in $\mathcal{K}$ is strongly Gaussian, and these criteria are
applied to find asymptotic estimations on the growth of the coefficients of
generating functions of combinatorial interest.

For admissible functions, Hayman proved (but the proof also works for strongly Gaussian functions)
the following central limit theorem (see \cite{Hayman} and \cite{CFFM1} for more details).

\begin{theor}[Hayman's local central limit theorem]
\label{teor:local central limit}
If $f(z)=\sum_{n=0}^\infty a_n z^n $ in $\K$ is strongly Gaussian, then
\begin{equation}
\label{eq:local central limit}
\lim_{t \uparrow R} \Big(\sup\limits_{n \in \Z} \Big|\frac{a_n t^n}{f(t)} \,\sqrt{2\pi}\,\sigma_f(t)-e^{-(n-m_f(t))^2/(2\sigma_f^2(t))}\Big|\Big)=0,
\end{equation}
where for $n<0$ it is understood that $a_n=0$.
\end{theor}

As a corollary of this theorem, we obtain the following.

\begin{coro}
\label{cor:m va a infinito si L1 gaussiana}
If $f \in \mathcal{K}$ is strongly Gaussian, then $f$ is a clan.
\end{coro}

\begin{proof} Restricting the supremum of Theorem \ref{teor:local central limit} to $n=-1$, we get
\[
\lim_{t\uparrow R} \exp\Big(-\frac{(m_f(t)+1)^2}{2\sigma_f(t)^2}\Big)=0.
\]
Since $\lim_{t \uparrow R}\sigma_f(t)=+\infty$, because $f$ is strongly Gaussian, we deduce that $f$ is a clan.
\end{proof}

This corollary shows that a large collection of functions ranging from $f(z)=e^z$ to the generating functions of
partitions $\prod_{j\ge 1}(1-z^j)^{-1}$, $|z|<1$, are clans; see \cite{CFFM1} and \cite{CFFM2} for these and other
interesting examples.

Notice that there are clans in $\K$ which are not strongly Gaussian. For instance, if $g\in\mathcal{K}$ is strongly
Gaussian, then $f(z)=g(z^N)$, with $N\in\mathbb N$, is a clan, but it is not strongly Gaussian since $a_{n}\equiv 0$
when $n$ is not a multiple of $N$, and hence \eqref{eq: asintotico coeficientes} cannot hold.

\subsubsection{Some power series that are not clans}
\label{section:examples not clans}

The geometric and negative binomial families \emph{are not clans}, since for $f(z)=1/(1-z)^N$, see Proposition \ref{prop:moments geometric},
we have that $\lim_{t \uparrow 1} \sigma_f(t)/m_f(t)=1/\sqrt{N}$.

In fact, for each $\alpha >0$, the function $f(z)=1/(1-z)^\alpha$, which is in $\K$,
is not a clan, because $\lim_{t \uparrow 1} \sigma_f(t)/m_f(t)=1/\sqrt{\alpha}$.

Even further, for {the function} $f(z)=1+\ln  (1/(1{-}z) )$, which is in $\K$, it holds that $\lim_{t \uparrow 1} \sigma_f(t)/m_f(t)=+\infty$.

Many other examples of power series in $\K$ which are not clans are provided by the following
immediate corollary of Theorem~\ref{toer:limsup del quotient}.
\begin{coro}
Let $f$ be in $\K$. If\/ $M_f=+\infty$ and $\overline{G}(f)>1$, then $f$ is \emph{not} a clan.
\end{coro}

For instance, power series with radius of convergence $R=1$, like $1+ \sum_{k=1}^\infty z^{2^k}$,
or entire power series like $1+\sum_{k=1}^\infty z^{2^k}/2^k!$, are in $\K$ but they are not clans.

Observe that for $f(z)=1/(1-z)$, which is not a clan, we have $M_f=+\infty$ and $\overline{G}(f)=1$.

\subsection{Some basic properties of clans}\label{subsection:properties clans}

We register now a few properties of power series in $\K$ that are clans, that is, functions in $\mathcal{K}$,
with radius of convergence $R>0$, and satisfying the limit condition given in \eqref{eq:lim para clan}.

\smallskip

$\bullet$ It follows from Chebyshev's inequality that if $(X_t)_{t\in[0,R)}$ is a clan, then for any $\varepsilon >0$,
\[
\lim_{t \uparrow R} \P\Big(\Big|\frac{X_t}{\E(X_t)}-1\Big|>\varepsilon\Big)=0,
\]
and thus that $X_t/\E(X_t)$ converges in probability to the constant 1 as $t\uparrow R$:
the random variable $X_t$ concentrates about its mean $m_f(t)$ as $t \uparrow R$.

\smallskip

$\bullet$ By Lemma \ref{lema:from integer moments to general moments}, we have that if $f$ is a clan, then
\[
\lim_{t \uparrow R} \frac{\E(X_t^p)}{\E(X_t)^p}=1\quad\text{for any $p \in (0,2]$}.
\]
Theorem \ref{teor:asintotico de momentos si clan} below will show that if $f$ is a clan,
then this limit result actually holds for any $p>0$.

\smallskip
$\bullet$ If $f$ is a clan (with at least three nonzero coefficients), then $\mathcal{D}_f=zf^\prime(z)$
is also a clan. This will be proved right after Theorem~\ref{teor:asintotico de momentos si clan}.
(The condition of three nonzero coefficients excludes the case in which the variables associated
to $\mathcal{D}_f$ are constant.)

$\bullet$ If $g$ is a clan, then for any integer $N\ge 1$, $f(z)=g(z^N)$ is also a clan,
since for $t\in[0,R^{1/N})$, $m_f(t)=Nm_g(t^N)$ and $\sigma_f(t)=N\sigma_g(t^N)$, where $m_g$, $\sigma^2_g$ and $m_f$,
$\sigma^2_f$ denote the mean and variance functions of the Khinchin families of $g$ and $f$.

\smallskip

$\bullet$ If $f$ and $g$ are clans with the same radius of convergence $R>0$, then their product $h\equiv fg$ is also a clan.
For we have $m_h=m_f + m_g$ and $\sigma_h^2=\sigma_f^2 + \sigma_g^2=o(m_f^2 + m_g^2)$ as $t \uparrow R$,
and thus $\sigma_h=o(m_h)$ as $t \uparrow R$. In particular, if $f$ is a clan
and $N\ge 1$ is an integer, then $f^N$ is also a clan.

\smallskip

$\bullet$ Finally, if $f$ and $g$ are entire functions which are clans,
then the composition $f\circ g$ is a clan. This is clear if both $f$ and $g$ are polynomials.
Otherwise, this follows from the identity
\[
L_{f\circ g}(t)=L_f(g(t))+\frac{L_g(t)}{m_f(g(t))},\quad \text{for any $t >0$},
\]
and by combining Lemma \ref{lem:condition for clan if Mfinfinite} and the fact that,
for a polynomial $h$ in $\K$ of degree $N$, we have that $\lim_{t \to \infty} L_h(t)=1-{1}/{N}$,
while $\lim_{t \to \infty} m_h(t)=N$.

In particular, if $g$ is an entire function which is a clan, then $e^g$ is a clan.

\subsection{Weak clans}\label{sec:weak clans}
 We say that $f \in \K$ is a \emph{weak clan} if
\[
\liminf_{t \uparrow R} \frac{ \E(X_t^2)}{\E(X_t)^2}=1, \quad \text{or, equivalently, if}
\quad \liminf_{t\uparrow R}\dfrac{\sigma_f(t)}{m_f(t)}=0.
\]
Of course, clans are weak clans.

\begin{prop}\label{prop:entire weak clans}
Every \emph{entire} function $f$ in $\K$ is a weak clan.
\end{prop}

As a consequence of this result, for entire functions $f \in \K$ which are not clans,
the quotient $\sigma_f(t)/m_f(t)$ must oscillate and has not limit as $t \to \infty$.

In contrast, power series $f \in \K$ with finite radius of convergence \emph{need not be weak clans}.
{Recall, for instance, the examples $f(z)=1/(1-z)$ and $f(z)=1+\ln (1/(1-z))$, for which $\lim_{t \uparrow 1} \sigma_f(t)/m_f(t)$ is 1 and $+\infty$, respectively.}


\begin{proof}[Proof of Proposition \textup{\ref{prop:entire weak clans}}]Polynomials are clans. For a transcendental entire function $f \in \K$ we have, because of Lemma \ref{lem:polynomial clans},  that $M_f=+\infty$. That $f$ is a weak clan follows then from combining \eqref{eq:liming L=1, transcendental} with \eqref{eq:clan in terms of f previo} of Lemma \ref{lema:previo weak clan}.\end{proof}

\begin{remark}  In fact, more is true. The following argument is based on  Rosenbloom's proof in~\cite{Rosenbloom} of the Wiman--Valiron theorem.

For a transcendental entire function $f$, we are going to show that, for every $\eta>0$ and $\varepsilon >0$, there exists a
set $G_\varepsilon \subset (0,R)$ of logarithmic measure
not exceeding $\varepsilon$, such that
\begin{equation}
\label{eq:limit weak}
\lim_{\substack{t \to \infty \\ t \notin G_\varepsilon}} \frac{\sigma^2_f(t)}{m_f(t)^{1+\eta}}=0,
\end{equation}
which, in particular, implies that
\[
\liminf_{t \to \infty} \frac{\sigma^2_f(t)}{m_f(t)^{1+\eta}}=0.
\]

For $a>0$, consider $H_a=\{x \ge a: \sigma^2_f(x)\ge m_f(x)^{1+\eta/2}\}$. For $x \in H_a$, we have that
\[
\frac{m_f^\prime(x)}{m_f(x)^{1+\eta/2}} \ge \frac{1}{x},
\]
and thus, using that $M_f=+\infty$, we deduce that
\[
\int_{H_a} \frac{dx}{x} \le \frac{2}{\eta \, m_f(a)^{\eta/2}}\,\cdot
\]
Let $a=a(\varepsilon)$ be such that $2/m_f(a)^{\eta/2}\le \eta \varepsilon$. Then, for $G_\varepsilon=H_{a(\varepsilon)}$,
we have that $G_\varepsilon$ has logarithmic measure at most $\varepsilon$, and for $t \notin G_\varepsilon$, we conclude  that
\[
\sigma^2_f(t) \le m_f(t)^{1+\eta/2}.
\]
Since $M_f=+\infty$, this gives \eqref{eq:limit weak}.
\end{remark}

\subsection{On the mean function of a clan}

Hayman showed in Lemmas 2 and 3 of~\cite{Hayman}) that the mean $m_f(t)$ of functions $f\in\mathcal{K}$
in the Hayman class cannot grow too slowly. Recall, from Section \ref{section:Strongly Gaussian Khinchin families},
that the Hayman class is a subclass of the class of strongly Gaussian functions, and thus, by
Corollary~\ref{cor:m va a infinito si L1 gaussiana}, functions in the Hayman class are clans.
Next, building upon Hayman's approach, we present a characterization of clans in terms of $m_f$ alone.

\smallskip

Observe that to first order approximation we have that
\[
m_f(t+t/m_f(t))\approx m_f(t)+m_f^\prime(t) \,\frac{t}{m_f(t)}=m_f(t)+\frac{\sigma_f^2(t)}{m_f(t)} ,
\]
and thus that
\[
\frac{m_f (t+{t}/{m_f(t)} )}{m_f(t)}\approx 1+\frac{\sigma_f^2(t)}{m_f(t)^2}\cdot
\]
This suggests that clans, for which 
$\lim_{t \uparrow R} \sigma_f(t)/m_f(t)=0$, may be characterized in terms of
the behaviour of the quotient $m_f (t+{t}/{m_f(t)} )/m_f(t)$ as $t \uparrow R$.
This is the content of Theorem \ref{teor:char de clans}.

Also,
 to second order approximation we have, using \eqref{eq:m and sigma in terms of f}, that
\begin{equation}
\label{eq:second order approximation}
\begin{aligned}
\ln f(t+t/m_f(t))-\ln f(t)
&\approx (\ln f)^\prime(t) \frac{t}{m_f(t)}+\frac{1}{2}\,(\ln f)^{\prime\prime}(t)\,
\frac{t^2}{m_f(t)^2}
\\
&=1+\frac{1}{2} \frac{\sigma_f^2(t)}{m_f(t)^2}-\frac{1}{2}\frac{1}{m_f(t)},
\end{aligned}
\end{equation}
which in turn suggests that clans, at least when $M_f=+\infty$, may be characterized in terms
of the behaviour of the difference $\ln f(t+{t}/{m_f(t)} )-\ln f(t)$, as $t \uparrow R$.
This is the content of Theorem \ref{teor:f crece lenta si f clan}.

\begin{theor}
\label{teor:char de clans}
Let $f\in\mathcal{K}$ have radius of convergence $0<R\le \infty$.

If\/ $R=\infty$, the power series $f$ is a clan if and only if
\begin{equation}
\label{eq:char de clan con m}
\lim_{t\uparrow R}\dfrac{m_f (t+{t}/{m_f(t)} )}{m_f(t)}=1.
\end{equation}

If\/ $R<\infty$, the power series $f$ is a clan if and only if \eqref{eq:char de clan con m}
holds, and besides,
\begin{equation}
\label{eq:char de clans con m con R finito}
\lim_{t\uparrow R}\,(R-t)\,m_f(t)=\infty.
\end{equation}
\end{theor}

\begin{remark}\label{remark:Borels lemma}
The classical Borel lemma, see for instance Chapter 9 of~\cite{Rubel}, claims that for any function $\mu(t)$
continuous and increasing in $[T, \infty)$ for some $T$, and if $a>1$, there exists an exceptional set $E$ of
logarithmic measure at most $a/(a - 1)$ so that
\[
\mu (t + t/\mu(t) )\le a \mu(t), \quad \text{for any $t \in [T, +\infty)\setminus E$}.
\]
\end{remark}

\begin{proof}[Proof of Theorem \textup{\ref{teor:char de clans}}]
If $f$ is a polynomial, then $\lim_{t \to \infty} m_f(t)=\text{deg}(f)$, and thus~\eqref{eq:char de clan con m} holds.
We may assume thus that $f$ is not a polynomial.

\smallskip

For the \emph{direct part}, we assume that $f$ is a clan and not a polynomial, and thus that $M_f=+\infty$.
Consider $t \in (0,R)$, and denote by $\Delta(t)$ the supremum
\[
\Delta(t)=\sup_{s\in[t,R)} \frac{\sigma_f^2(s)}{m_f(s)^2}\,\cdot
\]
Since $f$ is a clan, $\lim_{t \uparrow R} \Delta(t)=0$. Now, take $0<t \le r \le s<R$.
Since $\sigma_f^2(r)=rm_f'(r)$, we have that
\[
\frac{m_f^\prime(r)}{m_f(r)^2}\le \frac{\Delta(t)}{r}\,\cdot
\]
After integration in the interval $(t,s)$ and using that $\ln y \le y{-}1$, for any $y>1$, we get that
\begin{equation}
\label{eq:clans and means 1}
\frac{1}{m_f(t)}-\frac{1}{m_f(s)}\le \Delta(t) \ln \frac{s}{t}\le \frac{\Delta(t)}{t} \,(s-t).
\end{equation}
If $R<\infty$, by taking the limit as $s \uparrow R$, and using that $M_f=+\infty$, we get that
\begin{equation}
\label{eq:clans and means 2}
\frac{1}{m_f(t)} \le \frac{\Delta(t)}{t} (R-t), \quad \text{for $t \in [0,R)$}\,
\end{equation}
and since $\lim_{t \uparrow R} \Delta(t)=0$, we deduce that \eqref{eq:char de clans con m con R finito} holds.

If $R$ is finite, from \eqref{eq:clans and means 2} and because $\lim_{t\uparrow R} \Delta(t)=0$,
there exists $T\in(0,R)$ such that $t+t/m_f(t)<R$ for every $t\in(T,R)$; this, of course, is also true if $R=+\infty$.
In~\eqref{eq:clans and means 1}, take $t\in(T,R)$ and $s=t+t/m_f(t)<R$ and multiply by $m_f(t)>0$ in both sides to get
\[
0 \le 1 -\frac{m_f(t)}{m_f\big(t+{t}/{m_f(t)}\big)}\le \Delta(t).
\]
Then \eqref{eq:char de clan con m} follows since $\lim_{t \uparrow R}\Delta(t)=0$.

\smallskip

Next, the \emph{converse part} of the statement.

Notice first that we have that $M_f=+\infty$. For $R=+\infty$, this follows since $f$ is not a polynomial.
For $R<+\infty$,
this follows from the assumption \eqref{eq:char de clans con m con R finito}.

Observe that, it does not matter whether $R$ is finite or not,
we always have that $t+t/m_f(t)<R$ for $t\in [T,R)$, for appropriate $T\in (0,R)$, and
thus that $m_f(t+t/m_f(t))$ is well defined.

Denote $\lambda(t)=1+1/m_f(t)$. Consider the power series $\mathcal{D}_f$.
From Lemma \ref{lema:simic1} applied to~$\mathcal{D}_f$, we have that
\begin{equation}
\label{eq:clans and means 3}
m_{\mathcal{D}_f}(t) \ln \lambda(t) <\ln \frac{{\mathcal{D}_f}\big( \lambda(t) t\big)}{{\mathcal{D}_f}(t)},
\quad \text{for $t \in [T,R)$}.
\end{equation}

Since ${\mathcal{D}_f}(t)=m_f(t) f(t)$, for $t \in [0,R)$, appealing again to Lemma \ref{lema:simic1},
but now to $f$ itself, we have that, for $t \in [T,R)$,
\begin{equation}
\label{eq:clans and means 4}
\ln \frac{{\mathcal{D}_f}\big(\lambda(t) t\big)}{{\mathcal{D}_f}(t)}=\ln \frac{m_f\big(\lambda(t) t\big)}{m_f(t)}
+\ln \frac{f\big(\lambda(t) t\big)}{f(t)}\le \ln \frac{m_f\big(\lambda(t) t\big)}{m_f(t)} +m_f\big(\lambda(t) t\big) \ln \lambda(t).
\end{equation}
Combining inequalities \eqref{eq:clans and means 3} and \eqref{eq:clans and means 4}, dividing by $m_f(t) \ln \lambda(t)$,
and using equation~\eqref{eq:quotient mDF by mf}, we obtain that
\[
1\le \frac{\E(X_t^2)}{\E(X_t)^2}=\frac{m_{\mathcal{D}_f}(t)}{m_f(t)}
\le \frac{\ln \frac{m_f(\lambda(t) t)}{m_f(t)}}{m_f(t) \ln \lambda(t)} +\frac{m_f(\lambda(t) t)}{m_f(t)},
\quad \text{for $t\in[T,R)$}.
\]

Since $\lim_{t \uparrow R} m_f(\lambda(t) t)/m_f(t){=}1$, by the hypothesis \eqref{eq:char de clan con m},
and $\lim_{t \uparrow R} m_f(t) \ln \lambda(t){=}1$, because $\lim_{t \uparrow R} m_f(t)=\infty$, we deduce that
$\lim\nolimits_{t \uparrow R} \E(X_t^2)/\E(X_t)^2=1$, and thus that $f$ is a clan.
\end{proof}

\subsection{On the quotient \texorpdfstring{$f(\lambda t)/f(t)$}{f(lambdat)/f(t)}}

Let $f$ be in $\K$, with radius of convergence $R\le +\infty$. For $\lambda >1$, we consider the quotient $f(\lambda t)/f(t)$.
Bounds for this quotient appeared already in Lemma \ref{lema:simic1} of Simi{\'c}.

\begin{lemma}
\label{lema:f(lambda t)/f(t) increasing}
If $f$ is a power series in $\K$ with radius of convergence $R\le +\infty$,
then for $\lambda >1$, the function $f(\lambda t)/f(t)$ is increasing for $t\in (0,R/\lambda)$.
\end{lemma}

\begin{proof}
Fix $\lambda >1$. Consider the function $g(t)=\ln f(\lambda t)-\ln f(t)$, for $t\in (0,R/\lambda)$.
The derivative of $g$ multiplied by $t$ is $m_f(\lambda t)-m_f(t)$, which is always positive,
since $m_f$ is strictly increasing in $(0,R)$.
\end{proof}

\begin{lemma}
\label{lema:f(lambda t)/f(t) for a poly}
Let $f$ be an entire power series in $\K$ and let $\lambda >1$. If $f$ is a polynomial of degree $d$,
then $f(\lambda t)/f(t)$ is bounded for $t \in (0,\infty)$\textup{;}
in fact, $\lim_{t\to \infty}f(\lambda t)/f(t)=\lambda^{d}$. If~$\,f$ is transcendental,
then $\lim_{t\to \infty}f(\lambda t)/f(t)=+\infty$.
\end{lemma}

\begin{proof}
By Lemma \ref{lema:f(lambda t)/f(t) increasing}, we have that $\lim_{t\to \infty}f(\lambda t)/f(t)=\sup_{t>0}f(\lambda t)/f(t):=H$,
which could be $+\infty$.

If $H<+\infty$, we have, in particular, that $f(\lambda^n)\le H^n f(1)$, for each $n \ge 1$ and thus
that $f$ is a polynomial of degree at most $\ln H/\ln \lambda$. For if $a_k$ is the $k$-th coefficient
of $f$, then, by Cauchy estimates of coefficients, we have that $0\le a_k \le f(\lambda^n)/\lambda^{n\,k}\le f(1)H^n/\lambda^{n\,k}$,
for any $n \ge 1$, and so $a_k=0$ for any $k >\ln H/\ln \lambda$.
Conversely, if $f\in \K$ is a polynomial of degree $d$, then $\lim_{t \to \infty}f(\lambda t)/f(t)=\lambda^d$.
\end{proof}

Regarding Lemma \ref{lema:f(lambda t)/f(t) for a poly}, see item 24 in P{\' o}lya--Szeg\H{o} \cite{PolyaSzego},
and compare with Lemma~\ref{lema:simic1}.

\begin{lemma}
Let $f$ be an entire power series in $\K$. Assume that for a continuous function $\lambda(t)$ defined in $(0,+\infty)$
with values in $(1,+\infty)$, we have that $f(\lambda(t) t)/f(t)$ is bounded for $t \in (0,+\infty)$. Then
\[
\lambda(t)=1+O(1/m_f(t)), \quad \text{as $t \to \infty$}.
\]
\end{lemma}

\begin{proof}
Assume first that $f$ is not a polynomial, and thus that $M_f=+\infty$.
Let $J >0$ be such that $f(\lambda(t) t)/f(t)\le e^J$, for every $t>0$. Lemma \ref{lema:simic1} gives us that
\begin{equation}
\label{eq:bound m times ln lambda}
m_f(t) \ln \lambda(t) \le J, \quad \text{for $t >0$}.
\end{equation}
Define $\delta(t):= \lambda(t)-1$. Because of \eqref{eq:bound m times ln lambda}, we have that
\[
m_f(t) \delta(t)\le \big(e^{J/m_f(t)}-1\big) m_f(t)\, ,\quad \text{for $t >0$.}
\]
And so
\[
\limsup_{t \to \infty}m_f(t) \delta(t)\le \limsup_{t \to \infty} \big(e^{J/m_f(t)}-1\big) m_f(t)=J.
\]
The equality on the right holds because $M_f=+\infty$.

\smallskip

For a polynomial $f$, if $f(\lambda(t) t)/f(t)$ is bounded for $t \in (0,+\infty)$, then it follows immediately that $\lambda(t)$ must be  bounded. This is consistent with the conclusion of the lemma, since $\lim_{t \to \infty} m_f(t)=\text{degree}(f)$.
\end{proof}

The discussion above leads us to consider most naturally the case $\lambda(t)=1+1/m_f(t)$.

\begin{lemma}
\label{lemma:formula quotient flambda}
Let $f\in \mathcal{K}$ with radius of convergence $R\le\infty$. Assume that
\begin{equation}
\label{eq:t(1+1/m) less R}
\mbox{there exists $T\in (0,R)$ such that $t+t/m_f(t)<R$, for any $t\in[T,R)$.}
\end{equation}
Then
\begin{equation}
\label{eq:formula quotient flambda uno}
\frac{f(t+t/m_f(t))}{f(t)}=\sum_{k=0}^\infty \frac{1}{k!} \,\,\frac{\E(X_t^{\underline{k}})}{\E(X_t)^k},
\quad \text{for any $t \in (T,R)$}
\end{equation}
and, in particular,
\begin{equation}
\label{eq:formula quotient flambda dos}
\frac{f(t+t/m_f(t))}{f(t)}\ge \frac{1}{k!} \,\,\frac{\E(X_t^{\underline{k}})}{\E(X_t)^k} ,
\quad \text{for any $k \ge 0$.}
\end{equation}
\end{lemma}

The condition \eqref{eq:t(1+1/m) less R} in Lemma \ref{lemma:formula quotient flambda} is satisfied whenever $f$ is a clan.
To see this, observe first that if $R=+\infty$, then \eqref{eq:t(1+1/m) less R} is obvious. Now, if $R<+\infty$, and $f$ is a clan,
then \eqref{eq:char de clans con m con R finito} of Theorem \ref{teor:char de clans} gives, in particular, that $(R-t) m_f(t)>R>t$,
for any $t\in [T,R)$, for some $T\in(0,R)$, and thus that $t+t/m_f(t)<R$, for $t \in[T,R)$, which is~\eqref{eq:t(1+1/m) less R}.

\begin{proof}[Proof of Lemma \textup{\ref{lemma:formula quotient flambda}}]
Fix $t \in [T,R)$. The radius of convergence of the Taylor expansion of $f$ around $t$ exceeds $t/m_f(t)$, and this gives that
\begin{equation}
\label{eq:identity series normalized moments}
f(t+t /m_f(t))=\sum_{k=0}^\infty \frac{1}{k!}\,\,\frac{t^k f^{(k)}(t)}{m_f(t)^k}\,.
\end{equation}
Dividing by $f(t)$ and appealing to the expression \eqref{eq:formula factorial moments} of
the factorial moments of the $(X_t)$ in terms of $f$,
we obtain~\eqref{eq:formula quotient flambda uno}.
The inequalities in \eqref{eq:formula quotient flambda dos} follow since all the summands
in \eqref{eq:formula quotient flambda uno} are nonnegative.
\end{proof}

In the next result, as suggested by the second order approximation \eqref{eq:second order approximation},
clans are characterized as in Theorem \ref{teor:char de clans}, but involving $\ln f(t)$ instead of $m_f(t)$.

\begin{theor}\label{teor:f crece lenta si f clan}
Let $f\in \mathcal{K}$ with radius of convergence $R\le\infty$. Assume that $M_f=+\infty$. Then
$f$ is a clan if and only if condition \eqref{eq:t(1+1/m) less R} holds and
\begin{equation}
\label{eq:f crece lenta si f clan}
\lim_{t \uparrow R} \, \,\ln \Big(\frac{f (t+t/m_f(t) )}{f(t)}\Big)= 1.
\end{equation}
\end{theor}

Concerning the hypothesis $M_f=+\infty$ of Theorem \ref{teor:f crece lenta si f clan},
observe that for any polynomial $f\in \mathcal{K}$ of degree $N$, we have that
\[
\lim_{t \to \infty}\frac{f(t+t/m_f(t))}{f(t)}=\Big(1+\frac{1}{N}\Big)^N<e.
\]

\begin{proof} Assume first that $f$ is a clan.

We have observed  after the statement of Lemma \ref{lemma:formula quotient flambda} that condition \eqref{eq:t(1+1/m) less R} is satisfied by clans. To
verify \eqref{eq:f crece lenta si f clan}, denote $\lambda(t):=1+1/m_f(t)$, for $t \in (0,R)$.
Observe that since $\lim_{t \uparrow R} m_f(t)=M_f=+\infty$, we have that
\begin{equation}
\label{eq:lim m times ln lambda}
\lim_{t \uparrow R} m_f(t)\ln\lambda(t)=1.
\end{equation}
From Lemma \ref{lema:simic1}, we have, for $T \in (0, R)$ as in condition \eqref{eq:t(1+1/m) less R}, that
\[
m_f(t) \ln \lambda(t)\le \ln \frac{f(\lambda(t) t)}{f(t)}\le \frac{m_f(\lambda(t) t)}{m_f(t)} \, (m_f(t) \ln \lambda(t) ),
\quad \text{for any $t\in[T,R)$}.
\]
Using \eqref{eq:lim m times ln lambda} and \eqref{eq:char de clan con m} of Theorem \ref{teor:char de clans},
the limit \eqref{eq:f crece lenta si f clan} follows.

\smallskip

For the converse implication, assuming now that \eqref{eq:t(1+1/m) less R} \and \eqref{eq:f crece lenta si f clan} hold,
we will verify that $f$ is a clan. It is enough to show that
\[
\limsup_{t\uparrow R} \frac{\E(X_t^2)}{\E(X_t)^2}\le 1,
\]
 or, because of Corollary \ref{cor:growth of moment and factorial moment} and the hypothesis $M_f=+\infty$, that
\begin{equation}
\label{eq:clan con factorial dos} \limsup_{t\uparrow R} \frac{\E(X_t^{\underline{2}})}{\E(X_t)^2}\le 1.
\end{equation}

Using hypothesis \eqref{eq:t(1+1/m) less R} and Lemma \ref{lemma:formula quotient flambda} we obtain that
\[
\frac{f(t+t/m_f(t))}{f(t)}=\sum_{k=0}^\infty \frac{1}{k!} \frac{\E(X_t^{\underline{k}})}{\E(X_t)^k},
\quad \text{for $t \in [T,R)$.}
\]

Fix an integer $N\ge 3$. Since the summands above are all nonnegative we may bound
\[
\frac{f(t+t/m_f(t))}{f(t)}\ge \sum_{k=0}^N \frac{1}{k!} \frac{\E(X_t^{\underline{k}})}{\E(X_t)^k}
,\quad \text{for $t \in [T,R)$.}
\]
We split the sum on the right separating the summands corresponding to $k \le 2$ and those with $3\le k \le N$:
\[
\frac{f(t+t/m_f(t))}{f(t)}\ge 1+1+\frac{1}{2} \frac{\E(X_t^{\underline{2}})}{\E(X_t)^2}
+\sum_{k=3}^N \frac{1}{k!} \frac{\E(X_t^{\underline{k}})}{\E(X_t)^k}.
\]
For $3 \le k \le N$, we have that
\begin{equation}
\label{eq:liminf higher moments}
\liminf_{t \uparrow R } \frac{\E(X_t^{\underline{k}})}{\E(X_t)^k}
=\liminf_{t\uparrow R} \frac{\E(X_t^{\underline{k}})}{\E(X_t^k)} \,\frac{\E(X_t^k)}{\E(X_t)^k}\ge 1,
\end{equation}
since Jensen's inequality gives that $\E(X_t^k) \ge \E(X_t)^k$ and $\lim_{t \uparrow R}{\E(X_t^{\underline{k}})}/{\E(X_t^k)}=1$,
because of Corollary \ref{cor:growth of moment and factorial moment} and the hypothesis $M_f=+\infty$.

Fix now $\tau \in (0,1)$. Because of \eqref{eq:liminf higher moments}, there exists $S=S(N,\tau) \in [T,R)$ so
that $\E(X_t^{\underline{k}})/\E(X_t)^k\ge \tau$, for any $t \in [S,R)$ and any $3 \le k \le N$.
Thus, we have that
\[
\frac{f(t+t/m_f(t))}{f(t)}\ge 1+1+\frac{1}{2} \frac{\E(X_t^{\underline{2}})}{\E(X_t)^2}+\tau \sum_{k=3}^N \frac{1}{k!} ,
 \quad \text{for $t \in [S,R)$}.
\]
From this inequality and the hypothesis \eqref{eq:f crece lenta si f clan}, we deduce that
\[
e\ge 1+1+\frac{1}{2} \limsup_{t\uparrow R} \frac{\E(X_t^{\underline{2}})}{\E(X_t)^2}+\tau \sum_{k=3}^N \frac{1}{k!}\cdot
\]
Letting now $\tau \uparrow 1$, and then $N \to \infty$, we obtain that
\[
e\ge 1+1+\frac{1}{2} \limsup_{t\uparrow R} \frac{\E(X_t^{\underline{2}})}{\E(X_t)^2}+ \Big(e-1-1-\frac{1}{2}\Big),
\]
and conclude that \eqref{eq:clan con factorial dos} holds.
\end{proof}

In terms of the moment generation function of $X_t$, clans are characterized {as follows}. 

\begin{coro}\label{cor:f crece lenta si f clan bis}
Let $f\in \mathcal{K}$ with radius of convergence $R\le\infty$. Assume that $M_f=+\infty$. Then
$f$ is a clan if and only if condition \eqref{eq:t(1+1/m) less R} holds and
\begin{equation}
\label{eq:f crece lenta si f clan bis}
\lim_{t \uparrow R} \, \,\E\big(e^{(X_t-m_f(t))\nu(t)}\big)=1,
\end{equation}
where $\nu_f(t):=\ln(1+1/m_f(t))$ for $t>0$.
\end{coro}

\begin{proof} The results follows from the expression
\[
\E\big(e^{(X_t-m_f(t))\nu(t)}\big)=\frac{f(t+t/m_f(t))}{f(t)} \,\Big(1+\frac{1}{m_f(t)}\Big)^{-m_f(t)}
\]
and Theorem \ref{teor:f crece lenta si f clan}.
\end{proof}

Observe that, for polynomials in $\K$, \eqref{eq:f crece lenta si f clan bis} holds.

\subsection{Moments of clans}\label{subsection:moments of clans}
Our next result shows that for a clan $f$, any moment, not just the second one, is asymptotically equivalent,
as $t\uparrow R$, to the corresponding power of~$m_f(t)$.

\begin{theor}
\label{teor:asintotico de momentos si clan}
If $f\in \mathcal{K}$ with radius of convergence $R\le \infty$ is a clan with associated family $(X_t)_{t\in[0,R)}$, then
\begin{equation}
\label{eq:moment p behaves}
\lim_{t\uparrow R} \frac{\E(X_t^p)}{\E(X_t)^p}=1, \quad \text{for every $p>0$}.
\end{equation}
\end{theor}

As a consequence of Theorem \ref{teor:asintotico de momentos si clan}, we can now prove that, as anticipated
in Section~\ref{subsection:properties clans}, if~$f$ is a clan (with at least three nonzero coefficients),
then $\mathcal{D}_f=zf^\prime(z)$ is also a clan. Let $(X_t)_{t \in [0,R)}$ and $(W_t)_{t \in [0,R)}$ denote,
respectively, the Khinchin families of $f$ and of $\mathcal{D}_f$.

As observed in \eqref{eq:moments of derivative}, the moments of the families $(W_t)$ and $(X_t)$ are related by
\[
\E(W_t^p)=\frac{1}{m_f(t)} \,\E(X_t^{p + 1}), \quad \text{for any $p >0$ and any $t \in (0,R)$}.
\]
Thus,
\[
\frac{\E(W_t^2)}{\E(W_t)^2}= \frac{\E(X_t^3)}{m_f(t)}\, \frac{m_f(t)^2}{\E(X_t^2)^2} = \frac{\E(X_t^3)}{m_f(t)^3}\,
\Big(\frac{m_f(t)^2}{\E(X_t^2)}\Big)^2.
\]
Since $f$ is a clan, both fractions on the far right tend towards 1 as $t\uparrow R$ and,
consequently,~$\mathcal{D}_f$ is, as claimed, also a clan.

\smallskip

As another consequence, which we have also anticipated, 
observe that since the partition function $P(z)=\prod_{k=1}^\infty 1/(1-z^k)$ is a clan, Theorem \ref{teor:asintotico de momentos si clan}
and \eqref{eq:mean and variance partition} give for the moments of its associated family $(X_t)$
that, for any $p >0$,
\[
\E(X_t^p)\sim\E(X_t)^p\sim\frac{\zeta(2)^p}{(1-t)^{2p}} , \quad \text{as $t \uparrow 1$.}
\]

\begin{remark}\label{remark:momentos de Hayman}
Hayman, in Theorem III of~\cite{Hayman}, shows that the successive derivatives of Hayman (admissible) functions satisfy asymptotic
formulas which are equivalent to the conclusion of Theorem~\ref{teor:asintotico de momentos si clan}, $\lim_{t\uparrow R} \E(X_t^k)/\E(X_t)^k=1$
for $k\ge 1$ integer. Our probabilistic proof below shows that this conclusion is valid under the simple and more general notion of clan.
\end{remark}

\begin{proof}[Proof of Theorem \textup{\ref{teor:asintotico de momentos si clan}}]
Due to Lemma \ref{lema:from integer moments to general moments}, it is enough to prove \eqref{eq:moment p behaves} for any integer~$k~\ge~1$.

If $f$ is a polynomial of degree $N$, we have
that $\lim_{t \to \infty} \E(X_t^k)=N^k$ for any integer $k\ge 1$. In particular,
\[
\lim_{t\to \infty} \frac{\E(X_t^k)}{\E(X_t)^k}=\frac{N^k}{N^k}=1.
\]

We assume now that $f$ is not a polynomial and, consequently, that $M_f=+\infty$.

We first check that
\begin{equation}
\label{eq:cotas O momentos clan}
\limsup_{t \uparrow R} \frac{\E(X_t^k)}{\E(X_t)^k }\le e k!.
\end{equation}
Corollary~\ref{cor:growth of moment and factorial moment}, using that $M_f=+\infty$,
gives that the inequality \eqref{eq:cotas O momentos clan} is equivalent to
\[
\limsup_{t \uparrow R} \frac{\E(X_t^{\underline{k}})}{\E(X_t)^k }\le e k!,
\]
which follows from \eqref{eq:formula quotient flambda dos} and Theorem \ref{teor:f crece lenta si f clan}.

\smallskip

Denote $V_t:=X_t/\E(X_t)$, for $t\in(0,R)$, so that, for any integer $k \ge 1$ and any $t \in (0,R)$
\[
\frac{\E(X_t^k)}{\E(X_t)^k}= \E(V_t^k).
\]
We aim to show that for a clan, $\lim_{t\uparrow R} \E(V_t^k)=1$ for any integer $k\ge 1$.
For $k=1$, we have that $\E(V_t)=1$, for any $t \in (0,R)$, and the case $k=2$ is just the definition of clan.

By \eqref{eq:cotas O momentos clan}, for any integer $k\ge 1$, the moments of $V_t$ satisfy that
\begin{equation}
\label{eq:acotacion de momentos}
\limsup_{t \uparrow R} \E(V_t^{2k})\le e \, (2k)!.
\end{equation}

Fix an integer $k\ge 3$. Consider a constant $\omega >0$ and apply the Jensen, Cauchy--Schwarz and Chebyshev inequalities:
\[
\begin{aligned}
1&= \E(V_t)^k\le \E(V_t^k)=\E\big(V_t^k \uno\nolimits_{\{|V_t-1|>\omega\}}\big)+\E\big(V_t^k \uno\nolimits_{\{|V_t-1|\le\omega\}}\big)
\\
 &\le\E(V_t^{2k})^{1/2}\, \P\big(|V_t-1|>\omega\big)^{1/2}+(1+\omega)^k
\le \E(V_t^{2k})^{1/2}\, \frac{\sigma_f(t)}{m_f(t)}\,\frac{1}{\omega}+(1+\omega)^k,
\end{aligned}
\]
where with $\uno_A$ we denote the indicator function of the event $A$. Since $f$ is a clan, $\lim_{t \uparrow R}\sigma_f(t)/m_f(t)=0$, and this and the bound in \eqref{eq:acotacion de momentos}
combine to  imply that
\[
1\le \limsup_{t \uparrow R} \E(V_t^k)\le (1+\omega)^k,
\]
for any $\omega>0$. Therefore, $\limsup_{t \uparrow R} \E(V_t^k)=1$. Now since $k\ge 1$,
we have that $\E(V_t^k)\ge \E(V_t)^k=1$ for any $t\in(0,R)$, we get, as desired, that
\[
\lim_{t \uparrow R} \E(V_t^k)=1\quad\text{for every integer $k\ge 1$.}\qedhere
\]
\end{proof}

\begin{remark}
Theorem \ref{teor:asintotico de momentos si clan} does not hold for general families or sequences of random variables,
i.e., if $(V_n)_{n \ge 1}$ is a sequence of nonnegative random variables, then the condition $\lim_{n \to \infty} \E(V_n^2)/\E(V_n)^2=1$
does not imply that $\lim_{n \to \infty} \E(V_n^p)/\E(V_n)^p=1$ for $p>2$ (although this would be the case for $p<2$, because of
Lemma \ref{lema:from integer moments to general moments}).

Indeed, for $n \ge 1$, define $V_n$ taking the value $\sqrt{n}/\ln (n + 1)$ with probability $1/(n + 1)$ and the
value $\frac{1}{n}\big((n + 1)-\sqrt{n}/\ln(n + 1)\big)$ with probability $n/(n + 1)$.

For this sequence of random variables, we have that $\E(V_n)=1$ for any $n \ge 1$ and $\lim_{n \to \infty} \E(V_n^p)=1$,
if $p \le 2$, but $\lim_{n \to \infty} \E(V_n^p)=+\infty$, if $p>2$.
\end{remark}

\section{Order of entire functions and power series distributions}\label{section:entire functions}

In this section, we deal with entire functions $f$ in $\K$. In particular, we shall be interested
in the relation between the order (of growth) of $f$ and the growth of the mean $m_f(t)=\E(X_t)$,
of the variance $\sigma^2_f(t)=\V(X_t)$, and of moments $\E(X_t^p)$, with $p>0$, of the associated
family of probability distributions.

Recall that the order $\rho(f)$ of an entire function $f$ in $\K$ is given by
\begin{equation}
\label{eq:formula order}
\rho(f):= \limsup_{t\to \infty} \frac{\ln \ln \max_{|z|= t}\{|f(z)|\}}{\ln t}=\limsup_{t\to \infty} \frac{\ln \ln f(t)}{\ln t},
\end{equation}
where we have used \eqref{eq:max in t} in the second expression.
On the other hand, for any entire function $f(z)=\sum_{n=0}^\infty a_n z^n$, Hadamard's formula (see Theorem 2.2.2 in \cite{Boas})
gives $\rho(f)$ in terms of the coefficients of $f$:
\begin{equation}
\label{eq:Hadamard formula order}
\rho(f)=\limsup_{n \to \infty} \frac{n \ln n}{\ln(1/|a_n|)}\cdot
\end{equation}

\subsection{The order of \texorpdfstring{$f$}{f} entire and the moments \texorpdfstring{$\E(X_t^p)$}{E(Xt\textasciicircum p)}}

We can express the order of an entire function $f\in \K$ in terms of the mean $m_f(t)$ and, in fact, of any moment~$\E(X_t^p)$, as follows.

\begin{theor}
\label{teor:orden y m}
Let $f\in \mathcal{K}$ be an entire function of order $\rho(f)\le +\infty$. Then
\begin{equation}
\label{eq:order and moments}
\limsup_{t\to \infty} \frac{\ln [\E(X_t^p)^{1/p} ]}{\ln t}=\rho(f)\,, \quad \text{for any $p\ge 1$.}
\end{equation}
\end{theor}

The case $p=1$ of Theorem \ref{teor:orden y m}, i.e.,
$
\limsup_{t\to \infty} \ln m_f(t)/\ln t=\rho(f)$, appears in P\'olya and Szeg{\H o}, see item 53 in p. 9 of~\cite{PolyaSzego}.

\begin{proof} We abbreviate and write
\[
\Lambda_p:=\limsup_{t\to \infty} \frac{\ln [\E(X_t^p)^{1/p} ]}{\ln t}\, ,\quad \text{for $p>0$.}
\]
Observe that $\Lambda_p\le \Lambda_q$ if $0<p\le q$, by Jensen's inequality.

\smallskip

(a) First we show that
\[ \Lambda_p\le \rho(f)\,, \quad \text{for any $p >0$.}
\]
Fix $p>0$. The inequality trivially holds if $\rho(f)=+\infty$, so we may assume $\rho(f)<+\infty$.
Let $\omega>\rho(f)$ and take $\tau=(\omega+\rho(f))/2$. If $f(z)=\sum_{n=0}^\infty a_n z^n$, for $z \in \C$,
then Hadamard's formula \eqref{eq:Hadamard formula order} gives $N=N_\tau>0$ such that
\[
a_n\le \frac{1}{n^{n/\tau}},\quad \text{if } n\ge N.
\]

For $t$ such that $t^\tau>N$, we have
\[
\begin{aligned}
f(t) \,\E(X_t^p) &=\sum_{n=1}^\infty n^p a_n \,t^n=\sum_{n\le t^\omega} n^p a_n \,t^n+\sum_{n> t^\omega} n^p a_n \,t^n
\\
 &\le t^{p\omega} f(t)+\sum_{n\ge 1} n^p \,\frac{1}{n^{n/\tau}} \,n^{n/\omega}= t^{p\omega} f(t)+C\,,\end{aligned}
\]
where $C=C(\omega, \rho(f), p)<+\infty$. Thus,
\[
\E(X_t^p) \le t^{p\omega}+C/f(t),\quad \text{if } t^\tau>N,
\]
and so $\Lambda_p\le \omega$, for every $\omega >\rho(f)$, which implies, as desired, that $\Lambda_p\le \rho(f)$ for $p>0$ fixed above.

\smallskip

(b) To finish the proof, it is enough to show that
$ \rho(f)\le \Lambda_1
$,
because $\Lambda_1\le \Lambda_p$ for $p\ge 1$, and this, combined with part (a), would give \eqref{eq:order and moments}.

We may assume that $\Lambda_1<+\infty$, since otherwise there is nothing to prove.
We observe first that, for any $\omega>\Lambda_1$, there exists $T=T_\omega$ such that
\[
m_f(t)=\E(X_t)\le t^\omega , \quad\text{for } t\ge T,
\]
which, recall \eqref{eq:m and sigma in terms of f}, can be written in terms of $f$ as
\[
\frac{f^\prime(t)}{f(t)}< t^{\omega-1}, \quad\text{for } t\ge T.
\]
Upon integration, the above inequality gives that
\[
\ln f(t)-\ln f(T)\le \frac{1}{\omega} \,(t^\omega-T^\omega ), \quad\text{if } t\ge T,
\]
which implies, by the very definition \eqref{eq:formula order} of order, that
$\rho(f)\le \omega$. From this, we deduce, as desired, that~$\rho(f)\le \Lambda_1$.
\end{proof}

\begin{remark}
We mention that, if $\rho(f)$ is finite, then \eqref{eq:order and moments} holds also for $p\in (0,1)$.

To see this, fix $p\in (0,1)$. We just need to show that $\rho(f)\le \Lambda_p$. We are going to
interpolate $\Lambda_1$ between $\Lambda_p$ and $\Lambda_2$.

Let $u =1/(2-p) \in (1/2, 1)$, so that
\[
1=pu+2(1-u).
\]
By H\"older's inequality, we have that
\[
\E(X_t)=\E(X_t^{pu}\, X_t^{2(1-u)})\le \E(X_t^p)^u \,\E(X_t^2)^{1-u}, \quad \text{for any $t>0$,}
\]
and thus,
\[
\Lambda_1\le pu \,\Lambda_p+2(1-u)\,\Lambda_2.
\]
By \eqref{eq:order and moments}, $\Lambda_1=\Lambda_2=\rho(f)<+\infty$, so
\[
\rho(f) \le pu \,\Lambda_p+2(1-u)\,\rho(f).
\]
This, substituting the value of $u$, gives that $\rho(f)\le \Lambda_p$.
\end{remark}

\subsection{Entire gaps series, order and clans}
\label{section:entire gap series}
Recall, from Definition \ref{def:clan}, that an entire function $f$ in $\mathcal{K}$ is a clan if
\[
\lim_{t \to\infty}\frac{\sigma_f(t)}{m_f(t)}=0.
\]
We shall now exhibit examples of entire functions in~$\K$ of \emph{any given order} $ \rho$, $0\le \rho \le \infty$,
\emph{which are not clans}. These (counter)examples will be used in forthcoming discussions.

Fix $0<\rho<+\infty$.

Consider the increasing sequence of integers given by $n_1=0$, $ n_2=1$ and $n_{k + 1}=k n_{k}$, for any $k \ge 2$,
and let
\[
f(z)=1+\sum_{k=2}^\infty \frac{1}{n_k^{n_k/\rho}} \,z^{n_k}.
\]
The function $f$ is entire and belongs to $\K$. Hadamard's formula \eqref{eq:Hadamard formula order}
gives that $\rho(f) =\rho$.

Recall the gap parameter $\overline{G}(f)$ given in \eqref{eq:def de G barra}, and observe that,
in this case, $\overline{G}(f)={\limsup_{k\to\infty} k}=+\infty$, so from Theorem~\ref{toer:limsup del quotient},
we deduce that
\begin{equation}
\label{eq:limsup sigma over m}
\limsup_{t \to \infty}\frac{\sigma_f(t)}{m_f(t)}\ge 1
\end{equation}
holds, and, in particular, that $f$ is not a clan.

With the same specification of the sequence $n_k$, the power series $g$ and $h$ in $\K$ given by
\[
g(z)=1+\sum_{k=2}^\infty \frac{1}{n_k^{n_k^2}} \,z^{n_k}
\quad\text{and}\quad
h(z)=1+\sum_{k=2}^\infty \frac{1}{n_k^{n_k/\sqrt{\ln n_k}} }\, z^{n_k}
\]
are entire, of respective orders $\rho(g)=0$ and $\rho(h)=+\infty$, and are such
that \eqref{eq:limsup sigma over m} holds, and so, in particular, they are not clans.

\smallskip

The examples above of entire power series which are not clans are based on the seminal examples of Borel~\cite{Borel},
see also~\cite{Whittaker} of Whittaker, of entire functions whose lower order does not coincide with the order.

\subsection{The order of \texorpdfstring{$f$}{f} entire and the quotient \texorpdfstring{$\sigma^2_f(t)/m_f(t)$}{sigmaf\textasciicircum2(t)/mf(t)}}
The next result compares the order of the entire function with the quotient $\sigma_f^2(t)/m_f(t)$ as $t \to \infty$.

\begin{prop}
\label{prop:orden y sigma cuad partod por m}
Let $f\in \mathcal{K}$ be an entire function of order $\rho(f)\le \infty$. Then
\begin{align*}
\liminf_{t \to \infty} \frac{\sigma_f^2(t)}{m_f(t)}
&\le \liminf_{t \to \infty} \frac{\ln m_f(t)}{\ln t}\le \liminf_{t \to \infty} \frac{\ln \ln f(t)}{\ln t}
\\
&\le \limsup_{t \to \infty} \frac{\ln \ln f(t)}{\ln t}= \limsup_{t \to \infty} \frac{\ln m_f(t)}{\ln t}=\rho(f)
\le \limsup_{t \to \infty} \frac{\sigma_f^2(t)}{m_f(t)}\cdot
\end{align*}
\end{prop}

The equality statements in the middle of the comparisons of Proposition~\ref{prop:orden y sigma cuad partod por m}
are the very definition of order and the case $p=1$ of Theorem~\ref{teor:orden y m}.
On the other hand, B\'aez-Duarte shows in Proposition~7.7 of~\cite{BaezDuarteOtro} that
\[
\liminf_{t \to \infty} \frac{\sigma_f^2(t)}{m_f(t)}\le \rho(f) \le \limsup_{t \to \infty} \frac{\sigma_f^2(t)}{m_f(t)}\cdot
\]

\begin{proof}[Proof of Proposition \textup{\ref{prop:orden y sigma cuad partod por m}}]
We will check first that
\begin{equation}\label{eq:inequalities of limsups} \limsup_{t \to \infty} \frac{\ln \ln f(t)}{\ln t}
\le \limsup_{t \to \infty} \frac{\ln m_f(t)}{\ln t} \le \limsup_{t \to \infty} \frac{\sigma_f^2(t)}{m_f(t)}\cdot
\end{equation}
Of course, we already know that the first two limsup coincide with the order $\rho(f)$.

For the inequality on the right of \eqref{eq:inequalities of limsups}, let us denote
\[
L=\limsup_{t \to \infty} \dfrac{\sigma_f^2(t)}{m_f(t)}
\]
and assume that $L<+\infty$, since otherwise there is nothing to prove.
Recall that
\[
\frac{\sigma_f^2(t)}{m_f(t)}=\frac{t m_f^{\prime}(t)}{m_f(t)}\cdot
\]
Take $\omega >L$. Then there exists $T>0$ such that
\[
\frac{t m_f^{\prime}(t)}{m_f(t)}\le \omega, \quad\text{for any }t \ge T,
\]
and, by integration, for $t>T$,
\[
\ln m_f(t) -\ln m_f(T)\le \omega (\ln t - \ln T ),
\]
and thus
\[
\limsup_{t \to \infty} \frac{\ln m_f(t)}{\ln t}\le \omega,
\quad \text{and consequently,}\quad \limsup_{t \to \infty} \frac{\ln m_f(t)}{\ln t}\le L.
\]
The proof of the inequality on the left of \eqref{eq:inequalities of limsups},
\[
\limsup_{t \to \infty} \frac{\ln \ln f(t)}{\ln t}\le \limsup_{t \to \infty} \frac{\ln m_f(t)}{\ln t},
\]
follows as above using that
\[
t(\ln f(t))'=m_f(t).
\]

The proof of the inequalities for the liminf,
\[
\liminf_{t \to \infty} \frac{\sigma_f^2(t)}{m_f(t)} \le \liminf_{t \to \infty} \frac{\ln m_f(t)}{\ln t}
\le \liminf_{t \to \infty} \frac{\ln \ln f(t)}{\ln t},
\]
is analogous.
\end{proof}

Concerning Proposition \ref{prop:orden y sigma cuad partod por m}, a few observations are in order.

\medskip

(a) The three liminf in the statement, in general, do not give the order of $f$,
since the third one is the lower order of $f$, which, for instance and again,
by Borel~\cite{Borel}, does not have to coincide with the order.

(b) As for the third limsup in the second line: as pointed out by B\'aez-Duarte (see~\cite{BaezDuarteOtro}, p.~100), in general,
the order $\rho(f)$ is not given by $\limsup_{t \to \infty} \sigma_f^2(t)/m_f(t)$, as proposed by Kosambi in
Lemma~4 of~\cite{Kosambi2}. The example of Ba\'ez-Duarte is the canonical product $f$ given by
\begin{equation}
 \label{eq:example BD}
 f(z)=\prod_{n=1}^\infty \Big(1+\frac{z^{n^2}}{n^{4n^2}}\Big),
\end{equation}
that is an entire function in $\K$. Borel's theorem (see, for instance, Theorem 2.6.5 in~\cite{Boas}) tells us that
the order $\rho(f)$ of any canonical product coincides with the exponent of convergence of its zeros. The zeros of $f$
have exponent of convergence is $3/4$, and thus $f$ has order  $\rho(f)=3/4$. From the case $p=1$ of Theorem~\ref{teor:orden y m},
we see that $f$ satisfies that, say, $m_f(t)\le C\, t^{7/8}$, for some $C>0$ and every $t \ge 1$.
As $f(n^4 e^{i\pi/n^2})=0$, we obtain from Lemma \ref{lema:ceros caracteristica} that
 \[
 \sigma_f^2(n^4)\ge \frac{1}{4}\, n^4, \quad \text{for any $n \ge 1$},
 \]
and thus that
\[
 \frac{\sigma_f^2(n^4)}{m_f(n^4)}\ge \frac{1}{4C} \, n^{1/2}, \quad \text{for any $n \ge 1$}.
\]
So, for this particular function $f$, we have that $\limsup_{t \to \infty} \sigma_f^2(t)/m_f(t)=+\infty$, but
$\liminf_{t \to \infty} \sigma_f^2(t)/m_f(t)\le \rho(f)=3/4$.

Alternatively, and more generally, recall that for any $\rho \in [0, +\infty]$, we have exhibited in
Section \ref{section:entire gap series} an entire transcendental function $f \in \K$, with order $\rho(f)=\rho$,
and such that $\limsup_{t \to \infty} \sigma_f(t)/m_f(t)\ge 1$, and thus such that $\limsup_{t \to \infty} \sigma^2_f(t)/m_f(t)=+\infty$.
Also, $\liminf_{t \to \infty} \sigma^2_f(t)/m_f(t)\le \rho(f)<+\infty$, because of Proposition \ref{prop:orden y sigma cuad partod por m}.

\smallskip

(c) On the other hand, as a consequence of Proposition \ref{prop:orden y sigma cuad partod por m}, for an entire function $f \in \K$,
we have that \emph{if the limit $\lim_{t \to \infty}\sigma^2_f(t)/m_f(t)$ exists} (including the possibility of being $\infty$),
then the order $\rho(f)$ of $f$ is precisely
\begin{equation}\label{eq:order sigmacuad m} \rho(f)=\lim_{t \to \infty}\frac{\sigma^2_f(t)}{m_f(t)}\cdot
\end{equation}

In B\'aez-Duarte's example \eqref{eq:example BD}, the limit of ${\sigma^2_f(t)}/{m_f(t)}$ as $t\to\infty$ does not exist,
and~\eqref{eq:order sigmacuad m} does not hold. In general, \eqref{eq:order sigmacuad m} does not hold for any entire function $f$ in $\mathcal{K}$ for which the lower
order does not coincide with the order.

For the class of entire functions in $\K$ of genus zero which we are to discuss in the next Section~\ref{section:genus 0},
the identity~\eqref{eq:order sigmacuad m} holds (see the comments after Proposition~\ref{prop: funcion beta}).

The identity~\eqref{eq:order sigmacuad m} also holds for the class of nonvanishing entire functions in $\K$ of finite order.
To see this, let $f(z)=e^{g(z)} \in \K$, where $g$ is entire (not necessarily in $\K$).
We may assume that $g(t)\in \R$ for $t >0$, and further that, for some $T>0$,  $g(t)>0$~for $t\ge T$.
Assume that $f$ has finite order. Hadamard's factorization theorem gives that~$g$~is~a polynomial, say of {degree} 
 $N$.  Thus $\rho(f)=N$, by \eqref{eq:formula order}.
If the leading coefficient of the polynomial~$g$ is $b$, then
\begin{align*}
m_f(t)&=tf'(t)/f(t)= tg^\prime(t) \sim b N \,t^{N}, 
 \quad\text{as $t \to \infty$,}
 \\
\sigma^2_f(t)&=tm_f'(t)=t g^{\prime}(t)+t^2g^{\prime\prime}(t)\sim b N^2 \,t^{N}\quad\text{as $t \to \infty$.}
\end{align*}
Therefore, \eqref{eq:order sigmacuad m} holds.

\subsection{Entire functions of genus 0 with negative zeros}\label{section:genus 0}

We consider again the entire transcendental functions $f$ in $\K$ of genus $0$ (and thus of order $\le 1$) whose
zeros are all real and negative which we have considered in Section \ref{section:mean variance canonical products}. Recall that if normalized so that $f(0)=1$ they are the  canonical products of the form
\begin{equation}
\label{eq: producto infinito}
f(z)=\prod_{j=1}^\infty \Big(1+\frac{z}{b_j}\Big),\quad\text{for } z\in \C,
\end{equation}
and $(b_j)_{j\ge 1}$ is an increasing sequence of positive real numbers such that  $\sum_{j=1}^\infty {1}/{b_j}<\infty$.
The zeros $-b_1,-b_2,-b_3,\dots$ of $f$ all lie on the negative real axis.

 We {keep} 
 the notation $N(t)$ {for} 
 the {counting function} of the zeros of $f$:
\[
N(t)=\#\{j\ge 1: b_j\le t\}, \quad \mbox{for $t>0$}\,,
\]
which is a non-decreasing function such that $N(t)\to\infty$ as $t\to \infty$.

Recall also that Borel's theorem  tells us that the order $\rho(f)$ of the canonical product coincides with the exponent of convergence
of its zeros, which, in turn (see Theorem 2.5.8 in~\cite{Boas}), is given by
\[
 \limsup_{t\to\infty}\,\frac{\ln N(t)}{\ln t}=\rho(f).
\]

The function $\ln f(t)$, for $t \in (0,+\infty)$, may be expressed in terms of the counting
function $N(t)$; concretely, by integration by parts, one obtains
\begin{equation}
\label{eq:formula valiron}
\ln f(t)=\sum_{j=1}^\infty \ln\Big(1+\frac{t}{b_j}\Big)=\int_0^\infty \dfrac{tN(x)}{x(x+t)}\,dx=\int_0^\infty \dfrac{N(ty)}{y(y+1)}\,dy.
\end{equation}
This representation \eqref{eq:formula valiron} of $\ln f(t)$ in terms of the counting function $N(t)$ comes from {Valiron}, see~\cite{Valiron}.

From precise asymptotic information of the counting function $N(t)$ of $f$, one may obtain asymptotic information of the
mean and variance function of the family associated to $f$.

\begin{prop}
\label{prop: funcion beta}
Let $f$ be an entire function of genus $0$ with only negative zeros, and given by \eqref{eq: producto infinito}.
Assume that for $\rho\in(0,1)$ we have
\begin{equation}
\label{eq:asymptotic n}
N(t)\sim C t^{\rho}\, \quad \text{as $t\to\infty$}.
\end{equation} Then
\[
\begin{array}{rl}
\textup{(a)}\quad &\ln f(t)\sim \dfrac{C\pi}{\sin (\pi\rho )} \,t^\rho,\quad \text{ as } t\to\infty,\\[0.5cm]
\textup{(b)} \quad &m_f(t)\sim \dfrac{C\pi\rho}{\sin (\pi\rho )} \,t^\rho,\quad \text{ as } t\to\infty,\\[0.5cm]
\textup{(c)} \quad &\sigma_f^2(t) \sim \dfrac{C\pi\rho^2}{\sin(\pi\rho)} \,t^\rho,\quad \text{ as } t\to\infty.
\end{array}
\]

Conversely, if \textup{(a)}, \textup{(b)} or \textup{(c)} holds, then \eqref{eq:asymptotic n} holds.
\end{prop}

It follows that for entire functions of genus 0 given by formula \eqref{eq: producto infinito},
and whenever $N(t)\sim C t^{\rho}$ as $t\to \infty$, with $\rho \in (0,1)$, then
\[
{\lim_{t\to\infty} \frac{\sigma_f^2(t)}{m_f(t)}=\rho,}
\]
as was pointed out (but left unproved) by B\'aez-Duarte in Proposition 7.9 of~\cite{BaezDuarteOtro},
and also that equality holds in Proposition~\ref{prop:orden y sigma cuad partod por m}. Furthermore, it
follows that the entire functions of Proposition~\ref{prop: funcion beta} are clans; but see
Proposition~\ref{prop:canonical products with negative zeros} for a more general statement.

\begin{proof}
That $\textup{(c)}\Rightarrow \textup{(b)} \Rightarrow \textup{(a)}$ follows immediately by integration,
since $m_f(t)=t (\ln f)^\prime(t)$ and $\sigma^2_f(t)=t m_f^\prime(t)$; recall the formulas \eqref{eq:m and sigma in terms of f}.

That \textup{(a)} implies \eqref{eq:asymptotic n} is a classical Tauberian theorem,
see Valiron~\cite{Valiron} and Titchmarsh~\cite{Titchmarsh1, Titchmarsh2}.

The proof of the direct part of Proposition \ref{prop: funcion beta} follows from the
representation~\eqref{eq:formula valiron} and the following standard identity for the
Euler Beta function:
\begin{equation}
\label{eq:funcion beta}
\int_0^\infty\dfrac{y^\eta}{(1+y)^2}\,dy=\operatorname{Beta}(1+\eta,1-\eta)=\dfrac{\pi\eta}{\sin(\pi\eta)},
\quad \text{for any $\eta\in[0,1)$}.
\end{equation}

The representation \eqref{eq:formula valiron} gives that
\[
\ln f(t)=t^\rho \int_0^\infty \dfrac{y^\rho }{y(y+1)}\, \dfrac{N(ty)}{(ty)^\rho}\,dy.
\]
Since $N(ty)/(ty)^\rho\to C$ as $t\to \infty$ and $y^\rho/(y(y+1))$ is integrable in $[0,\infty)$,
\[
\lim_{t\to\infty}\dfrac{\ln f(t)}{t^\rho}=C\int_0^\infty\dfrac{y^\rho}{y(y+1)}\,dy.
\]
Integrating by parts, using that $\rho<1$ and \eqref{eq:funcion beta}, we obtain that
\[
\int_0^\infty\dfrac{y^\rho}{y(y+1)}\,dy
=\dfrac{1}{\rho}\int_0^\infty\dfrac{y^\rho}{(y+1)^2}\,dy
=\dfrac{\pi}{\sin(\pi\rho)}\cdot
\]
That is,
\[
\ln f(t)\sim \dfrac{C\pi}{\sin(\pi\rho)}\,t^\rho,\quad \text{as } t\to\infty.
\]

For the mean and the variance, we have the representations
\[
m_f(t)=t^\rho\int_0^\infty \dfrac{y^\rho }{(y+1)^2}\, \dfrac{N(ty)}{(ty)^\rho}\,dy
\quad \text{and} \quad \sigma_f^2(t)
=t^\rho\int_0^\infty \dfrac{y^\rho(y-1) \, }{(y+1)^3}\,\dfrac{N(ty)}{(ty)^\rho}\,dy.
\]
Arguing as above, (b) and (c) follow from \eqref{eq:asymptotic n}.
\end{proof}

\begin{remark}
For general canonical products $f$ of genus $p\ge 0$ with only negative zeros, there is a representation
of $\ln f(t)$ analogous to \eqref{eq:formula valiron}, which is the case $p=0$ (see, for instance,
Theorem 7.2.1 in~\cite{Bingham}):
\[
\ln f(t)=(-1)^p \int_0^\infty \frac{(t/x)^{p + 1}}{1 + t/x} \,N(x)\,\frac{dx}{x}\cdot
\]
This expression means in particular that for $p$ odd, the canonical product $f$ is bounded by 1, for $t\in (0,\infty)$,
and thus shows that $f$ is not in $\K$. Observe, in any case, that for a primary factor $E_p(z)$, which for $|z|<1$
is $E_p(z)=\exp(-\sum_{j=p + 1}^\infty z^j/j)$, we have that its $(2p+2)$-coefficient is $-p/(2(p+1)^2)$, and
therefore $E_p(-z/a)$ with $a>0$, which vanishes at $-a$, is not in $\K$, if $p\ge 1$; in general, thus,
canonical products of nonzero genus with only negative zeros are not in $\K$.
\end{remark}

\begin{remark}
For entire functions $f$ of genus zero and only negative zeros, if the number of zeros is \emph{comparable}
to a power $t^\rho$ with $\rho \in (0,1)$ ($N(t)\asymp t^{\rho}$, as $t\to\infty$), so are the mean and variance
 of its Khinchin family. This follows most directly from the representation~\eqref{eq:formula valiron}.
\end{remark}

Entire functions of genus 0 with only negative zeros are always clans.

\begin{prop}
\label{prop:canonical products with negative zeros}
Every entire function $f$ in $\K$ defined by \eqref{eq: producto infinito}
with $\sum_{j\ge 1} 1/b_j<+\infty$ is a clan.
\end{prop}

\begin{proof}
Assume that $f$ is not a polynomial.
For $f$ given by \eqref{eq: producto infinito}, we have,
recalling \eqref{eq:variance less than mean of canonical products}, that
\[
\sigma_f^2(t)<m_f(t),
\]
and thus $\sigma_f^2(t)/m_f^2(t)\le 1/m_f(t)$.
Since $m_f(t)\to\infty$ as $t\to\infty$, we obtain that $f$ is a clan.
\end{proof}

\subsection{Entire functions, proximate orders and clans}

For $\rho \ge 0$, a \emph{proximate $\rho$-order} $\rho(t)$ is a continuously differentiable
function defined in $(0,+\infty)$ and such that
\[
\lim_{t \to\infty} \,\rho(t)=\rho \quad \text{and} \quad \lim_{t \to \infty} \,\rho^\prime(t) t \ln t=0.
\]
Traditionally, proximate orders are allowed to have a discrete set of points where they are not
differentiable, but have both one-sided derivatives at those points, see Section 7.4 in~\cite{Bingham}.

If for a proximate $\rho$-order $\rho(t)$ we write $V(t)=t^{\rho(t)}$, for $t>0$, then
(see, for instance, Lemma 5 in Section 12, Chapter I, of~\cite{Levin2}),
for every $\lambda >0$ we have that
\begin{equation}
\label{eq:behaviour of V}
\lim_{t \to \infty} \frac{V(\lambda t)}{V(t)}=\lambda^\rho.
\end{equation}
In other terms, the function $V(t)/t^\rho$ is a \emph{slowly varying} function.

\begin{theor}[Valiron's proximate theorem for $\K$]
\label{theo:Valiron}
If $f$ is an entire function in $\K$ of finite order $\rho\ge 0$,
then there is a proximate $\rho$-order $\rho(t)$ such that
\[
\limsup_{t \to \infty}\frac{\ln f(t)}{t^{\rho(t)}}=1.\]
\end{theor}

See Theorem 7.4.2 in~\cite{Bingham}; the smoothness which we require in our definition of
proximate $\rho$-order $\rho(t)$ is provided by Proposition 7.4.1 and Theorem 1.8.2
(the smooth variation theorem) in \cite{Bingham}.

We have the following.

\begin{theor}
\label{theor:proximate order and clans}
Let $f$ be an entire function in $\K$ of finite order $\rho >0$, let $\rho(t)$ be a proximate $\rho$-order, and let $\tau >0$. Then
\begin{equation}
\label{eq:proximate ln f}
\lim_{t \to \infty}\frac{\ln f(t)}{t^{\rho(t)}}=\tau
\end{equation}
if and only if
\begin{equation}
\label{eq:proximate mf}
\lim_{t \to \infty}\frac{m_f(t)}{t^{\rho(t)}}=\tau \rho.
\end{equation}

If either \eqref{eq:proximate ln f} or \eqref{eq:proximate mf} holds, then $f$ is a clan.
\end{theor}

Observe that in both \eqref{eq:proximate ln f} and \eqref{eq:proximate mf} a `$\lim$' is assumed,
and not a `$\limsup$' as in Valiron's Theorem~\ref{theo:Valiron}, which encompass all finite order
entire functions. Condition \eqref{eq:proximate ln f} of Theorem \ref{theor:proximate order and clans}
concerns entire functions that are said to have \emph{regular growth}.

Comparing with Valiron's theorem, the limit $\tau$ instead of $1$ amounts no extra generality,
since replacing $\rho(t)$ by $\rho^\star(t)=\rho(t)+\ln \tau/\ln t$, we have that $\rho^\star(t)$
is also a proximate $\rho$-order and
\[
\lim_{t \to \infty}\frac{\ln f(t)}{t^{\rho^\star(t)}}=1.\]

That \eqref{eq:proximate mf} implies that $f$ is a clan is due to Simi{\'{c}, \cite{Simic}.
In \cite{Simic}, see also \cite{Simic0}, it is claimed, in the terminology of the present paper,
that any entire function of \emph{finite order} in $\mathcal{K}$ is a clan, which is not the case;
see, for instance, Section~\ref{section:entire gap series}. The error in the argument originates
in a misprint in the statement of Theorem 2.3.11 in p.~81 of~\cite{Bingham}: the $\limsup$ appearing
in that statement should be a $\liminf$ (which is what is actually proved in \cite{Bingham}).
See also \cite{Liu}, p.~101, for a similar warning. The argument of \cite{Simic}
shows precisely that~\eqref{eq:proximate mf} implies that~$f$ is a clan.

As for the implication $\eqref{eq:proximate ln f} \Rightarrow \eqref{eq:proximate mf}$,
compare with Lemma 3.1 in~\cite{OstrovskiiUreyen2}.

For constant proximate $\rho$-order, ($\rho(t)=\rho$, for $t>0$) that \eqref{eq:proximate ln f}
implies \eqref{eq:proximate mf} and then that $f$ is a clan is due to P\'olya and Szeg\"o with
an argument involving some delicate estimates: combine items 70 and 71, of page 12, of theirs \cite{PolyaSzego}.

\begin{proof}
Fix $\lambda>1$. Lemma \ref{lema:simic1} gives us that
\[
 m_f(t) \ln \lambda \le\ln f(\lambda t)-\ln f(t)\le m_f(\lambda t) \ln \lambda, \quad\mbox{for $t>0$},
 \]
and thus dividing by $V(t)=t^{\rho(t)}$,
\begin{equation}
\label{eq:m divided by V}
\frac{m_f(t)}{V(t)} \le \frac{1}{\ln \lambda}\, \Big[\frac{\ln f(\lambda t)}{V(\lambda t)}\,\frac{V(\lambda t)}{V(t)} -\frac{\ln f(t)}{V(t)}\Big]
\le \frac{m_f(\lambda t)}{V(\lambda t)}\,\frac{V(\lambda t)}{V(t)} , \quad\mbox{for $t>0$}.
\end{equation}

We first prove that $\eqref{eq:proximate ln f} \Rightarrow \eqref{eq:proximate mf}$.
From~\eqref{eq:behaviour of V}, \eqref{eq:proximate ln f} and~\eqref{eq:m divided by V}, and letting $t\to \infty$,
we deduce that
\[
\limsup_{t \to \infty}\frac{m_f(t)}{t^{\rho(t)}}\le \tau\, \frac{\lambda^\rho-1}{\ln \lambda}
\le \lambda^\rho \,\liminf_{t \to \infty}\,\frac{m_f(t)}{t^{\rho(t)}}\cdot
\]
Letting $\lambda \downarrow1$, equation \eqref{eq:proximate mf} follows.

\smallskip
We now prove that $\eqref{eq:proximate mf} \Rightarrow \eqref{eq:proximate ln f}$.
Assume first that
\begin{equation}
\label{eq:ln f divided by V}
\limsup_{t \to \infty} \frac{\ln f(t)}{t^{\rho(t)}}<+\infty.
\end{equation}
If \eqref{eq:ln f divided by V} holds, then from the first inequality of~\eqref{eq:m divided by V},
and using \eqref{eq:proximate mf}, we deduce that
 \[
 \tau \rho \ln \lambda +\liminf_{t\to \infty}\frac{\ln f(t)}{t^{\rho(t)}}
 \le \lambda^\rho \liminf_{t\to \infty}\frac{\ln f(t)}{t^{\rho(t)}},
 \]
while the second inequality of~\eqref{eq:m divided by V} gives
\[
 \lambda^\rho \limsup_{t\to \infty}\frac{\ln f(t)}{t^{\rho(t)}}\le \limsup_{t\to \infty}\frac{\ln f(t)}{t^{\rho(t)}}
 + \tau \rho \lambda^\rho \ln \lambda.
\]
Writing the two inequalities above as
\[
\tau \rho \le \frac{\lambda^\rho-1}{\ln \lambda} \, \liminf_{t\to \infty}\frac{\ln f(t)}{t^{\rho(t)}}
\quad \text{and}\quad
\frac{\lambda^\rho-1}{\ln \lambda} \, \limsup_{t\to \infty}\frac{\ln f(t)}{t^{\rho(t)}}\le \tau \rho \lambda^\rho ,
\]
and by letting $\lambda\downarrow 1$, we get that
\[
  \tau \le \liminf_{t\to \infty}\frac{\ln f(t)}{t^{\rho(t)}}
  \quad\text{and}\quad
  \limsup_{t\to \infty}\frac{\ln f(t)}{t^{\rho(t)}}\le \tau,
\]
so \eqref{eq:proximate ln f} follows.

\smallskip

To show that \eqref{eq:ln f divided by V} holds, we first observe that
\[
\frac{d}{ds} \,s^{\rho(s)}= (s \ln s\cdot \rho^\prime(s)+\rho(s) ) \,s^{\rho(s)-1}, \quad \text{for $s>0$},
\]
and from the defining properties of the proximate orders, we deduce that,
for appropriately large $A>0$,
\begin{equation}
\label{eq:cota derivada proximate}
\frac{d}{ds} \,s^{\rho(s)}\ge \frac{\rho}{2} \,s^{\rho(s)-1}, \quad \text{for each $s\ge A$}.
\end{equation}

From \eqref{eq:proximate mf}, by incrementing $A$ if necessary, we deduce that
\[
\frac{d}{ds}\ln f(s)=\frac{m_f(s)}{s}\le 2 \tau \rho \, s^{\rho(s)-1}\le 4 \tau \,\frac{d}{ds} \,s^{\rho(s)}, \quad \text{for $s\ge A$},
\]
where \eqref{eq:cota derivada proximate} was used in the last inequality. We thus have that
\[
\ln f(t)\le \ln f(A) +4 \tau \,t^{\rho(t)}-4 \tau A^{\rho(A)}, \quad \text{for $t \ge A$},
\]
from which we obtain that
\[
\limsup_{t \to \infty} \frac{\ln f(t)}{t^{\rho(t)}}\le 4\tau,
\]
as wanted.

This argument shows in fact that
\[
\frac{1}{\rho}\,\liminf_{t \to \infty} \frac{m_f(t)}{t^{\rho(t)}}\le \liminf_{t \to \infty} \frac{\ln f(t)}{t^{\rho(t)}}
\le \limsup_{t \to \infty} \frac{\ln f(t)}{t^{\rho(t)}}\le\frac{1}{\rho}\,\limsup_{t \to \infty} \frac{m_f(t)}{t^{\rho(t)}}\cdot
\]
For another proof of this last chain of inequalities, see, for instance,
Theorem 4 in~\cite{VaishKasana}.

\smallskip
Finally, we prove that $\eqref{eq:proximate mf} \Rightarrow \mbox{$f$ is a clan}$.
From \eqref{eq:behaviour of V} and \eqref{eq:proximate mf}, and taking into account
that $\rho>0$, we deduce that
\begin{equation}
\label{eq:regularly varying}
\lim_{t\to \infty}\frac{m_f(\lambda t)}{m_f(t)}=\lambda^\rho, \quad \text{for any $\lambda >0$},
\end{equation}
and thus that the mean $m_f$ is a regularly varying function, see Section~1.4 of~\cite{Bingham}.

\smallskip

To show that $f$ is a clan, we may assume that $f$ is not a polynomial,
and thus that $\lim_{t \to \infty} m_f(t)=+\infty$.

Now, \eqref{eq:regularly varying} implies that for any function $\lambda(t)$ such that $\lambda(t) >1$
and such that $\lim_{t \to \infty} \lambda(t)=1$, we have that
\begin{equation}
\label{eq:regularly varying bis}
\lim_{t \to \infty} \frac{m_f(\lambda(t)\, t)}{m_f(t)}=1.
\end{equation}
To see this, fix $\varepsilon >0$. Then we have that $1<\lambda(t)\le 1+\varepsilon$,
for $t\ge t_\varepsilon$. Therefore,
\[
1\le \frac{m_f(\lambda(t)\, t)}{m_f(t)}\le \frac{m_f((1 + \varepsilon)\, t)}{m_f(t)},
\quad \text{for $t \ge t_\varepsilon$}.
\]
Thus, $\limsup_{t \to \infty}{m_f(\lambda(t)\, t)}/{m_f(t)}\le (1+\varepsilon)^\rho$,
and thus \eqref{eq:regularly varying bis} holds.

\smallskip
Applying \eqref{eq:regularly varying bis} with $\lambda(t)=1+{1}/{m_f(t)}$ and appealing to
Theorem \ref{teor:char de clans}, we conclude that $f$ is a clan.
\end{proof}

\subsection{Exceptional values and clans}\label{section:exceptional values}

The entire gap series of Section \ref{section:entire gap series}, which are our basic examples
of entire functions in $\K$ which are not clans, have no Borel exceptional values. This follows,
for instance, from a classical result of Pfluger and P\'olya~\cite{PflugerPolya}. We show next that,
in general, entire functions which are not clans have no Borel exceptional values.

Recall that, by definition, $a$ is a Borel exceptional value of an entire function $f$ of finite order
if the exponent of convergence of the $a$-values of $f$ (i.e., the zeros of $f(z)-a=0$) is strictly
smaller than the order of $f$; a theorem of Borel claims that a nonconstant entire function of finite
{order} can have at most one Borel exceptional value.

\begin{theor}
\label{teor:Borel exceptional and clans}
If the entire function $f\in \K$ has finite order and has one Borel exceptional value, then $f$ is a clan.
\end{theor}

\begin{proof} Let $\rho$ be the order of $f$. Assume that $a\in \C$ is the Borel exceptional value for $f$.
Denote with $s$ the exponent of convergence of the zeros of $f(z)-a$. Thus $s<\rho$, since $a$ is a Borel
exceptional value for $f$.

Let $f(z)=a+P(z)\,e^{Q(z)}$ be the Hadamard factorization of $f$, where $P$ is the canonical product
formed with the zeros of $f-a$, and $Q$ is a polynomial of degree $d$ and leading coefficient $c \neq 0$.
The order of $P$ is $s$, and thus the order $\rho$ of $f$ must be the integer $d$.

Now,
\[
|f(t)-a|=|P(t)|e^{\Re Q(t)} \quad \text{and} \quad \ln |f(t)-a|=\ln |P(t)|+\Re Q(t), \quad \text{for $t>0$}.
\]

Take $s^\prime \in (s,d)$. For a certain $t^\prime$ depending on $s^\prime$, we have
for $t\ge t^\prime$ that $\ln |P(z)|\le |z|^{s^\prime}$, if $|z|= t$. Besides,
\[
\frac{\Re Q(t)}{t^d}=\Re c+O\Big(\frac{1}{t}\Big), \quad \text{as $t \to \infty$}.
\]
We conclude that
\[
\lim_{t \to \infty} \frac{\ln |f(t)-a|}{t^d}=\Re c,
\]
and therefore that
\[
\lim_{t \to \infty} \frac{\ln f(t)}{t^d}=\Re c.
\]

Observe that if $\Re c=0$, then $\Re Q(t)=O(t^{d-1})$ as $t \to\infty$, and that would mean that
\[
\limsup_{t \to \infty} \frac{\ln f(t)}{t^h}=0,
\]
for some $h$ such that $(s<)\,h<d$, and thus, in particular, that $f$ would be of order at most $h$,
which is not the case.

Thus $\Re c >0$, and condition \eqref{eq:proximate ln f} of Theorem \ref{theor:proximate order and clans}
holds, so $f$ is a clan.
\end{proof}

For an entire function, not necessarily of finite order, a Picard exceptional value is a value that
is taken just a finite number of times. For Picard exceptional values and clans, we have the following
result, which came out in a conversation of one of the authors with Walter Bergweiler.
\begin{prop}
\label{prop:exponential of finite order}
If $f=Pe^g$ is in $\mathcal{K}$, where $P$ is a polynomial and $g$ is an entire function
in $\K$ {\upshape{of finite order}}, then $f$ is a clan.
\end{prop}

The value $0$ is Picard exceptional for $f=Pe^g$. It is not assumed that $P$ is in $\mathcal{K}$,
but it is assumed that $g$ is in $\mathcal{K}$. Observe also that the assumption is that $g$ is of
finite order; if $e^g$ were of finite order, that $f$ is a clan would follow from
Theorem \ref{teor:Borel exceptional and clans}.

\begin{proof} The entire function $f$ is transcendental, since $g\in \mathcal{K}$ is not a constant.
From Lemma \ref{lem:condition for clan if Mfinfinite}, we have that $f$ is a clan if
and only if $\lim_{t \to \infty} L_f(t)=1$. To show this, we verify first that $g$ satisfies
\begin{equation}
\label{eq:condition for clan if entire exp} \lim_{t \to \infty} \frac{g^{\prime\prime}(t)}{g^\prime(t)^2}=0.
\end{equation}

Condition \eqref{eq:condition for clan if entire exp} clearly holds if $g$ is a polynomial.

Assume thus that $g$ is not a polynomial. From the case $p=1$ of Theorem~\ref{teor:orden y m} applied to
the derivative~$g^\prime$, which is also of finite order, it follows that for some finite constant $S>0$
and radius $R_1>0$, we have that
\[\frac{g^{\prime\prime}(t)}{g^\prime(t)} \le t^S, \quad \text{for $t>R_1$}.
\]
Besides, since $g^\prime $ is not a polynomial,  we have, for some radius $R_2$, that $g^\prime(t)>t^{S+1}$,
for $t >R_2$. And thus \eqref{eq:condition for clan if entire exp} holds.

Next, a calculation, recall \eqref{eq:def de Lf}, gives that
\[
L_f(t)=\Big(\dfrac{P''(t)}{P(t)}\dfrac{1}{g'(t)^2}+2\dfrac{P'(t)}{P(t)}\dfrac{1}{g'(t)}+\dfrac{g''(t)}{g'(t)^2}+1\Big){\Big/}
\Big(\dfrac{P'(t)}{P(t)}\dfrac{1}{g'(t)}+1\Big)^2\,, \quad \text{for $t >0$.}
\]
Since $P$ is a polynomial, we have that $P'(t)/P(t)$ and $P''(t)/P(t)$ tend to 0 as $t\to \infty$.
Besides, since $g \in \K$, we have that $\lim_{t \to \infty} g'(t)=+\infty$.
Using now~\eqref{eq:condition for clan if entire exp}, it is deduced that $\lim_{t \to \infty} L_f(t)=1$.
\end{proof}

\section{Some questions}
\label{section:questions}

(1) If $f\in\mathcal{K}$ has radius of convergence $R\le \infty$ and its Khinchin
family $(X_t)_{t\in[0,R)}$ satisfies for \emph{some} value of $p\ge 2$ that
\begin{equation}
\label{eq:question limit quotients}\lim_{t\uparrow R} \frac{\E(X_t^p)}{\E(X_t)^p}=1,
\end{equation}
then $f$ is a clan. This follows directly from Lemma \ref{lema:from integer moments to general moments}.
Assume that \eqref{eq:question limit quotients} is satisfied for \emph{some} value $p \in (1,2)$.
Is this enough to deduce that $f$ is a clan?

\smallskip

(2) Is it the case that, for every entire function $f$ in $\K$,
\[
\limsup_{t \to \infty} \frac{\sigma_f(t)}{m_f(t)}\le 1?
\]

\smallskip

(3) If $g$ is an entire function in $\K$ (not necessarily of finite order),
it is natural to ask whether  $f=e^g$ is always a clan or not. This is the
case if $g$ is a clan or if $g$ has finite order, as we have seen in Section~\ref{subsection:properties clans}
and in Proposition \ref{prop:exponential of finite order}.

\

\noindent\textbf{Funding.} A. Cant\'on has been partially supported by the grant PID-2021-124195NB-C31 from the Agencia
Estatal de Investigaci\'on of the Ministerio de Ciencia, Innovaci\'on y Universidades of Spain. J.\,L. Fern\'andez and P. Fern\'andez
acknowledge generous support of Fundaci\'on Akusmatika. Research of V.\,J. Maci\'a was partially funded by grants MTM2017-85934-C3-2-P2
and PID2021-124195NB-C32 of Ministerio de Econom\'{\i}a y Competitividad of Spain, by the European Research Council Advanced Grant 834728,
by the Madrid Government (Comunidad de Madrid-Spain) under programme PRICIT, and by Grant UAM-Santander for the mobility of
young researchers 2023.


\begin{thebibliography}{99}





\bibitem{Abi}
  Abi-Khuzam, F.\,F.:
  \href{https://doi.org/10.2307/2047294}{Maximum modulus convexity and the location of zeros of an entire function}.
  \emph{Proc. Amer. Math. Soc.} \textbf{106} (1989), no.~4, 1063--1068.


\bibitem{Abi2}
  Abi-Khuzam, F.\,F.:
  \href{https://doi.org/10.2307/2154769}{Hadamard convexity and multiplicity and location of zeros}.
  \emph{Trans. Amer. Math. Soc.} \textbf{347} (1995), no.~8, 3043--3051.


\bibitem{BaezDuarteOtro}
  B\'aez-Duarte, L.:
  \href{https://doi.org/10.1006/aima.1995.1028}{General Tauberian theorems on the real line}.
  \emph{Adv. Math.} \textbf{112} (1995), no.~1, 56--119.


\bibitem{BaezDuarte}
  B\'aez-Duarte, L.:
  \href{https://doi.org/10.1006/aima.1997.1599}{Hardy--{R}amanujan's asymptotic formula for partitions and the central limit theorem}.
  \emph{Adv. Math.} \textbf{125} (1997), no.~1, 114--120.


\bibitem{Bingham}
  Bingham, N.\,H., Goldie, C.\,M. and Teugels, J.\,L.:
  \emph{\href{https://doi.org/10.1017/CBO9780511721434}{Regular variation}}.
  Encyclopedia of Mathematics and its Applications 27, Cambridge University Press, Cambridge, 1987.

\bibitem{Boas}
  Boas, R.\,P., Jr.:
  \emph{Entire functions}.
  Academic Press, New York, 1954.


\bibitem{BoichukGoldberg}
 {Bo\u{\i}chuk}, V.\,S. and Gol'dberg, A.\,A.:
  \href{https://doi.org/10.1007/bf01153540}{The three lines theorem}.
  \emph{Math. Notes} \textbf{15} (1974), 26--30.

\bibitem{Borel}
  Borel, E.:
  \href{https://doi.org/10.1007/bf03013527}{Sur quelques fonctions enti\`{e}res}.
  \emph{Rend. Circ. Matem. Palermo} \textbf{23} (1907), 320--323.


\bibitem{CFFM1}
  Cant\'on, A., Fern\'andez, J.\,L., Fern\'andez, P. and Maci\'a, V.\,J.:
  \href{https://doi.org/10.1007/s40315-021-00420-6}{Khinchin families and Hayman class}.
  \emph{Comput. Methods Funct. Theory} \textbf{21} (2021), no.~4, 851--904.


\bibitem{CFFM2}
  Cant\'on A., Fern\'andez J.\,L., Fern\'andez P. and Maci\'a V.:
  \href{https://doi.org/10.1007/s00009-023-02579-9}{Khinchin families, set constructions, partitions and exponentials}.
  \emph{Mediterr. J. Math.} \textbf{21} (2024), article no.~39, 28~pp.

\bibitem{Rubel}
  Colliander, J.\,E. and Rubel, L.\,A.:
  \emph{\href{https://doi.org/10.1007/978-1-4612-0735-1}{Entire and meromorphic functions}}.
  Universitext, Springer, New York, 1996.

\bibitem{Hayman}
  Hayman, W.\,K.:
  \href{https://doi.org/10.1515/crll.1956.196.67}{A generalisation of {S}tirling's formula}.
  \emph{J. Reine Angew. Math.} \textbf{196} (1956), 67--95.

\bibitem{Hayman2}
  Hayman, W.\,K.:
  \href{https://doi.org/10.1090/pspum/011/0252639}{Note on {H}adamard's convexity theorem}.
  In \emph{Entire functions and related parts of analysis ({P}roc. {S}ympos. {P}ure {M}ath., {L}a {J}olla, {C}alif., 1966)}, pp. 210--213.
  Proc. Sympos. Pure Math.~11, Amer. Math. Soc., Providence, RI, 1968.

\bibitem{Hilberdink}
  Hilberdink, T.:
  \href{https://doi.org/10.1093/qmathj/haz061}{Asymptotics of entire functions and a problem of {H}ayman}.
  \emph{Q.\,J. Math.} \textbf{71} (2020), no.~2, 667--676.

\bibitem{JohnsonUnivariate}
  Johnson, N.\,L., Kemp, A.\,W. and Kotz, S.:
  \emph{\href{https://doi.org/10.1002/0471715816}{Univariate discrete distributions}}.
  Third edition.
  Wiley Series in Probability and Statistics, Wiley-Interscience [John Wiley \& Sons], Hoboken, NJ, 2005.

\bibitem{Kjellberg}
  Kjellberg, B.:
  The convexity theorem of Hadamard--Hayman.
  In \emph{Proceedings of the Symposium in Mathematics, Royal Institute of Technology, Stockholm, (June 1973)}, pages 87--114. 1973.

\bibitem{Kosambi1}
  Kosambi, D.\,D.:
  Characteristic properties of series distributions.
  \emph{Proc. National Inst. Sci. India} \textbf{15} (1949), 109--113.

\bibitem{Kosambi2}
  Kosambi, D.\,D.:
  Classical Tauberian theorems.
  \emph{J. Indian Soc. Agric. Stat.} \textbf{10} (1958), 141--149.

\bibitem{Levin2}
  Levin, B.\,J.:
  \emph{\href{https://bookstore.ams.org/mmono-5}{Distribution of zeros of entire functions}}.
  Revised edition.
  Translations of Mathematical Monographs 5, American Mathematical Society, Providence, RI, 1980.

\bibitem{Liu}
  Liu, Q.:
  \href{https://doi.org/10.1239/aap/1035227993}{Fixed points of a generalized smoothing transformation and applications to the branching random walk}.
  \emph{Adv. in Appl. Probab.} \textbf{30} (1998), no.~1, 85--112.

\bibitem{Noak}
  Noak, A.:
  \href{https://doi.org/10.1214/aoms/1177729894}{A class of random variable with discrete distribution}.
  \emph{Ann. Math. Statist.} \textbf{21} (1950), 127--132.

\bibitem{OstrovskiiUreyen2}
  Ostrovskii, I.\,V. and \"{U}reyen, A.\,E.:
  \href{https://doi.org/10.1080/0278107031000120431}{Distance between a maximum modulus point and zero set of an entire function}.
  \emph{Complex Var. Theory Appl.} \textbf{48} (2003), no.~7, 583--598.

\bibitem{OstrovskiiUreyen}
  Ostrovski{\u{\i}}, I.\,V. and \"{U}reyen, A.\,E.:
  \href{https://doi.org/10.1007/s10688-006-0047-7}{The growth irregularity of slowly growing entire functions}.
  \emph{Funct. Anal. Appl.} \textbf{40} (2006), 304--312.

\bibitem{Pitman}
  Pitman, J.:
  \href{https://doi.org/10.2307/2974785}{Some probabilistic aspects of set partitions}.
  \emph{Amer. Math. Monthly} \textbf{104} (1997), no.~3, 201--209.

\bibitem{PflugerPolya}
  Pfluger, A. and P\'olya, G.:
  \href{https://doi.org/10.1017/s0305004100013244}{On the power series of an integral function having an exceptional value}.
  \emph{Proc. Cambridge Phil. Soc.} \textbf{31} (1935), 153--155.

\bibitem{PolyaSzego}
  P\'olya, G. and Szeg{\H{o}}, G.:
  \emph{Problems and theorems in analysis. Vol. {II}}.
  Springer Study Edition, Springer, New York-Heidelberg, 1976.

\bibitem{Rosenbloom}
  Rosenbloom, P.\,C.:
  Probability and entire functions.
  In \emph{Studies in mathematical analysis and related topics}, pp. 325--332.
  Stanford Studies in Mathematics and Statistics~IV, Stanford Univ. Press, Stanford, CA, 1962.

\bibitem{Rossberg}
  Rossberg, H.-J.:
  \href{https://doi.org/10.1007/BF02363232}{Positive definite probability densities and probability distributions}.
  \emph{J. Math. Sci.} \textbf{76} (1995), no.~1, 2181--2197.

\bibitem{Sakovic}
  Sakovi{\v{c}}, G.\,N.:
  Pro {\v{s}}irinu spectra (On the width of the spectrum).
  \emph{Dopov{\={\i}}d{\={\i}} Akad. Nauk. Ukra\"{\i}n.} \textbf{11} (1965), 1427--1430.

  \bibitem{Schumitzky} Schumitzky, A.: \href{https://doi.org/10.1080/17476938908814380}{A probabilistic approach to the Wiman--Valiron theory for functions of several complex entire functions}. \emph{Complex Variables, Theory and Application.} \textbf{13} (1989), 85--98.

\bibitem{Simic}
  Simi{\'c}, S.:
  \href{https://doi.org/10.2174/1876527000901010003}{On moments of the power series distributions}.
  \emph{Open Stat. Prob. J.} \textbf{1} (2009), 3--6.

\bibitem{Simic0}
  Simi{\'c}, S.:
  \href{https://doi.org/10.1134/S0001434607050148}{Some properties of entire functions with nonnegative {T}aylor coefficients}.
  \emph{Math. Notes} \textbf{81} (2007), 681--685.

\bibitem{Titchmarsh1}
  Titchmarsh, E.\,C.:
  \href{https://doi.org/10.1112/jlms/s1-1.1.35}{A theorem on infinite products}.
  \emph{J. London Math. Soc.} \textbf{1} (1926), no.~1, 35--37.

\bibitem{Titchmarsh2}
  Titchmarsh, E.\,C.:
  \href{https://doi.org/10.1112/plms/s2-26.1.185}{On integral functions with real negative zeros}.
  \emph{Proc. London Math. Soc. (2)} \textbf{26} (1927), 185--200.

\bibitem{VaishKasana}
  Vaish, S.\,K. and Kasana, H.\,S.:
  On the proximate type of an entire function.
  \emph{Publ. Inst. Math. (Beograd) (N.S.)} \textbf{32(46)} (1982), 167--174.

\bibitem{Valiron}
  Valiron, G.:
  \href{https://doi.org/10.5802/afst.287}{Sur les fonctions enti{\`{e}}res d'ordre nul et d'ordre fini et en particulier les fonctions {\`{a}} correspondance r\'eguli{\`{e}}re}.
  \emph{Ann. Fac. Sci. Toulouse Sci. Math. Sci. Phys. (3)} \textbf{5} (1913), 117--257.

\bibitem{Whittaker}
  Whittaker, J.\,M.:
  \href{https://doi.org/10.1112/jlms/s1-8.1.20}{The lower order of integral functions}.
  \emph{J. London Math. Soc.} \textbf{8} (1933), no.~1, 20--27.


\end{thebibliography}
\end{document}